\documentclass[smallextended,numbook,runningheads]{svjour3}     
\smartqed  
\usepackage{graphicx}
\usepackage{amsmath,amsfonts,fancyhdr}
\usepackage{amssymb}
\usepackage{epstopdf} 
\usepackage{mathptmx}      
\usepackage{geometry}
\usepackage{latexsym}
\usepackage{booktabs,arydshln,multirow}
\usepackage{setspace}
\usepackage{xcolor}
\usepackage{stmaryrd}
\usepackage{cite}
\journalname{**********}
\begin{document}

\title{Sharp error analysis for averaging Crank-Nicolson schemes with corrections for subdiffusion with nonsmooth solutions
\thanks{
* Corresponding author.
}
}

\titlerunning{Sharp error analysis for ACN schemes}        

\author{Baoli Yin$^1$         \and
        Yang Liu$^1$  \and
        Hong Li$^{1,*}$ 
}


\institute{
Baoli Yin
\\
\email{baolimath@126.com}
\\
\\
Yang Liu
\\
\email{mathliuyang@imu.edu.cn}
\\
\\
Hong Li
\\
\email{smslh@imu.edu.cn}
\\
\\
${}^1$~ School of Mathematical Sciences, Inner Mongolia University, Hohhot 010021, China
}

\date{Received: date / Accepted: date}

\maketitle

\begin{abstract}
Thanks to the singularity of the solution of linear subdiffusion problems, most time-stepping methods on uniform meshes can result in $O(\tau)$ accuracy where $\tau$ denotes the time step.
The present work aims to discover the reason why some type of Crank-Nicolson schemes (the averaging Crank-Nicolson scheme) for the subdiffusion can only yield $O(\tau^\alpha)$$(\alpha<1)$ accuracy, which is much lower than the desired.
The existing well developed error analysis for the subdiffusion, which has been successfully applied to many time-stepping methods such as the fractional BDF-$p (1\leq p\leq 6)$, all requires singular points be out of the path of contour integrals involved.
The averaging Crank-Nicolson scheme in this work is quite natural but fails to meet this requirement. 
By resorting to the residue theorem, some novel sharp error analysis is developed in this study, upon which correction methods are further designed to obtain the optimal $O(\tau^2)$ accuracy.
All results are verified by numerical tests.
\end{abstract}

\section{Introduction}\label{intro}
We study numerically solving the subdiffusion problem on a polygon $\Omega \subset \mathbb{R}^d$ modeled by
\begin{equation}\label{intr.0}
  \partial_t^\alpha u(t)-\Delta u=f(x,t),\quad
  u(x,0)=u_0,\quad
  u(x,t)=0 \text{ for } x \in \partial\Omega,
  0<t\leq T,
\end{equation}
where $\alpha \in (0,1)$ denotes the order of the Caputo fractional derivative $\partial_t^\alpha$, which is defined by
\begin{equation*}
\partial_t^\alpha \phi(t)
=\frac{1}{\Gamma(1-\alpha)}
\int_{0}^{t}\frac{\phi'(s)}{(t-s)^\alpha}
\mathrm{d}s.
\end{equation*}
The equation such as (\ref{intr.0}) depicts the so called anomalous diffusion, marked by sublinear growth in mean-squared particle displacements.
Over the past few decades, it has gained considerable prominence for its broad utility in physics, biology, material science, finance, etc.; see \cite{barkai2000continuous,kilbas2006theory,magin2004fractional,golding2006physical,raberto2002waiting,zaslavsky2002chaos} and references cited therein.

Addressing the singularity issue in fractional models with Caputo fractional derivatives poses a formidable challenge when designing high-order numerical methods.
It is a well known result that the solution of subdiffusion is singular at initial time even though the initial condition and source terms are smooth \cite{stynes2016too,sakamoto2011initial}.
Several different techniques have been proposed to surmount such obstacle in literature.
The L-type methods such as the L1 method\cite{sun2006fully,lin2007finite}, L2-$1_\sigma$ method \cite{alikhanov2015new} and so on which were initially presented on uniform meshes and suffer the order-reduction problem, were extended to graded/nonuniform meshes, or modified by adding correction terms (still on uniform meshes). 
The authors in \cite{stynes2017error} carried out the analysis for finite difference method on graded meshes for subdiffusion with initial singularity.
Liao et al. in \cite{liao2018sharp} explored the L1 method on nonuniform meshes for reaction-subdiffusion equations by resorting to the discrete fractional Gr\"onwall inequality and obtained the sharp error estimates.
See also \cite{liao2021second} where a second-order difference method on nonuniform meshes was developed.
A framework for the analysis of the error of L1-type discretizations on graded temporal meshes in the $L^\infty$ and $L^2$ norm was presented by Kopteva in \cite{kopteva2019error} where both the finite difference methods and finite element methods were considered.
By carefully choosing the steps for nonuniform meshes, Chen and Styne in \cite{chen2019error}  studied the Alikhanov's high-order scheme for initial value/initial-boundary value problems with second-order accuracy.
Yan et al. in \cite{yan2018analysis} proposed the modified L1 scheme on uniform meshes for subdiffusion with smooth and nonsmooth initial data and derived the optimal accuracy.
See also the correction method for high-order L-$k$ scheme developed in \cite{shi2022correction}.
It is also notable that Li et al. in \cite{li2021novel} developed two type difference methods on uniform meshes by the variable changing technique and derived the error analysis for time-fractional parabolic equations.
By using the log orthogonal functions, Chen et al. in \cite{chen2020spectrally} developed the spectral-Galerkin method for solving the subdiffusion with spectrally accuracy.
\par
Another group of difference methods in fractional calculus discretization stems from the work \cite{lubich1986discretized} known as the convolution quadrature (CQ).
Since the CQ is developed on uniform meshes and in general correction terms are needs to maintain the high-order accuracy of the method.
In \cite{zeng2015numerical}, Zeng et al. considered two second-order fractional linear multistep methods (belonging to the CQ framework) with correction terms at each time step for the subdiffusion and derived the uniform optimal accuracy.
The authors in \cite{jin2017correction} developed correction formulas for the $k$th-order fractional backward difference formulas by adding correction terms at the initial $k-1$ steps and obtained the optimal accuracy at fixed positive time.
\par
Thanks to the initial singularity of the solution, most time-stepping methods such as the L1 \cite{yan2018analysis} or CQ \cite{jin2017correction} methods, can merely result in first-order accuracy in temporal direction at any positive time (see the convergence orders in the last column in Table \ref{tab1}).
However, a much lower accuracy is observed in the previous work \cite{yin2022efficient} (see Example 4 in \cite{yin2022efficient} with parameter $\theta$ there replaced by $\frac{1}{2}$), where the following Crank-Nicolson scheme is considered,
\begin{equation}\begin{split}\label{intr.1}
D_\tau^\alpha(U_h^{n-\frac{1}{2}}-u_h(0))
-\frac{1}{2}(\Delta_h U_h^n
+\Delta_h U_h^{n-1})
=\frac{1}{2}(f_h^n+f_h^{n-1}).
\end{split}\end{equation}
The term
$D_\tau^\alpha(U_h^{n-\frac{1}{2}}-u_h(0))$, known as a discrete convolution approximating the Caputo derivative $\partial_t^\alpha u$ at time $t_{n-\frac{1}{2}}$, is formulated by the natural averaging technique
\begin{equation}\begin{split}\label{intr.2}
D_\tau^\alpha(U_h^{n-\frac{1}{2}}-u_h(0))=
\frac{1}{2}
\big[
\mathcal{D}_\tau^\alpha (U_h^n-u_h(0))
+\mathcal{D}_\tau^\alpha (U_h^{n-1}-u_h(0))
\big],
\end{split}\end{equation}
in which $\mathcal{D}_\tau^\alpha (U_h^n-u_h(0))$ stands for any second-order difference formulas for $\partial_t^\alpha u$ at time $t_{n}$, defined by
\begin{equation*}\begin{split}
\mathcal{D}_\tau^\alpha (U_h^n-u_h(0))=
\tau^{-\alpha}
\sum_{k=0}^n \widetilde{\omega}_{n-k}(U_h^k-u_h(0)),\quad \text{and}\quad
\widetilde{\omega}(\zeta)=\sum_{k=0}^\infty \widetilde{\omega}_k \zeta^k.
\end{split}\end{equation*}
For example, one can take $\widetilde{\omega}(\zeta):=(\frac{3}{2}-2\zeta+\frac{1}{2}\zeta^2)^\alpha$ (the second-order fractional backward difference formula or the FBDF-2) or $\widetilde{\omega}(\zeta):=(1-\zeta)^\alpha[1+\frac{\alpha}{2}(1-\zeta)]$ (the generalized 2nd-order Newton–Gregory formula or the GNG-2) or other CQ methods which can be found in \cite{lubich1986discretized}.
In this work, we call the Crank-Nicolson scheme (\ref{intr.1}) combined with the averaging technique (\ref{intr.2}) the \textbf{averaging Crank-Nicolson (ACN) scheme}, as appeared in the title.
Despite possessing several advantages such as conciseness and symmetry, the ACN scheme, however, suffers from a significant drawback in terms of precision, as illustrated in Table \ref{tab1}, where only $O(\tau^\alpha)$ accuracy (see the convergence orders of ACN by FBDF-2 or GNG-2) can be observed even through the initial condition and source terms are smooth (in this case, $f=0$, $u_0(x)=\sin x$ and $u(x,t)$ is taken as $E_\alpha(-t^\alpha)\sin x$ with $E_\alpha(\cdot)$ denoted as the Mittag-Leffler function \cite{li2015numerical}).
To the best of our knowledge, the underlying causes leading to this low convergence accuracy remain presently elusive, which motivates us to develop sharp error analysis for such numerical scheme.
\begin{table}[]
\caption{The accuracy of the ACN scheme and FBDF-2 at $t=0.5$.}
\label{tab1}
\begin{tabular}{cccccccccc}
\hline\noalign{\smallskip}
\multirow{2}{*}{$\alpha$} & \multirow{2}{*}{$\tau$} & \multicolumn{2}{c}{ACN(FBDF-2)} &  & \multicolumn{2}{c}{ACN(GNG-2)} &  & \multicolumn{2}{c}{FBDF-2} \\ \cline{3-4} \cline{6-7} \cline{9-10} 
                          &                         & $\|U_h^n-u(t_n)\|_{L^2(\Omega)}$             & Order       &  & $\|U_h^n-u(t_n)\|_{L^2(\Omega)}$            & Order      &  & $\|U_h^n-u(t_n)\|_{L^2(\Omega)}$          & Order     \\ \hline \noalign{\smallskip}
\multirow{3}{*}{0.1}      & $2^{-7}$  & 4.3729E-01        &  --           &  & 4.2992E-01       &  --          &  & 2.4643E-04     &  --         \\
                          & $2^{-8}$  & 4.1777E-01        & 0.07        &  & 4.1056E-01       & 0.07       &  & 1.2305E-04     & 1.00      \\
                          & $2^{-9}$  & 3.9869E-01        & 0.07        &  & 3.9166E-01       & 0.07       &  & 6.1480E-05     & 1.00      \\ \hline \noalign{\smallskip}
\multirow{3}{*}{0.5}      & $2^{-7}$  & 5.3041E-02        &  --           &  & 5.0128E-02       &  --          &  & 1.3513E-03     & --          \\
                          & $2^{-8}$  & 3.7978E-02        & 0.48        &  & 3.5868E-02       & 0.48       &  & 6.7405E-04     & 1.00      \\
                          & $2^{-9}$  & 2.7096E-02        & 0.49        &  & 2.5578E-02       & 0.49       &  & 3.3663E-04     & 1.00      \\ \hline \noalign{\smallskip}
\multirow{3}{*}{0.9}      & $2^{-7}$  & 4.5622E-03        &  --           &  & 4.4804E-03       &  --          &  & 2.7063E-03     & --          \\
                          & $2^{-8}$  & 2.4458E-03        & 0.90        &  & 2.4020E-03       & 0.90       &  & 1.3470E-03     & 1.01      \\
                          & $2^{-9}$  & 1.3111E-03        & 0.90        &  & 1.2876E-03       & 0.90       &  & 6.7198E-04     & 1.00      \\ \hline
\end{tabular}
\end{table}
\par
In history, other types of Crank-Nicolson scheme have been studied.
In \cite{gao2015stability}, the authors considered the so-called \textbf{fractional Crank-Nicolson scheme} for the subdiffusion by approximating the derivative $\partial^\alpha_t u$ at point $t_{n-\frac{\alpha}{2}}$ by using the shifted fractional Euler method, and derived optimal second-order accuracy for sufficiently smooth solutions.
It is notable that the shifted fractional Euler method was first proposed by Dimitrov \cite{dimitrov2013numerical} and was
generalized within the shifted CQ framework \cite{liu2021unified}.
Considering the nonsmoothness of the solution, Jin et al. \cite{jin2018analysis} analyzed the fractional Crank-Nicolson scheme in detail for the subdiffusion $(\alpha <1)$ and proposed two-step correction methods to maintain the optimal accuracy.
Later on, Wang et al. \cite{wang2021single} improved this method by adding only one correction term at the first time step and still obtained the optimal accuracy.
Recently, the authors \cite{yin2022efficient} carried out rigorous error analysis for some type of corrected difference $\theta$-schemes by using formulas introduced in \cite{yin2020necessity}, which approximate the fractional derivative at $t_{n-\theta}$ with $\theta<\frac{1}{2}$.
We emphasize that although the corrected $\theta$-scheme \cite{yin2022efficient} can preserve the optimal second-order accuracy for $\theta<\frac{1}{2}$, the result and analysis can not be directly extended to the case $\theta=\frac{1}{2}$, mainly due to singularities of some key functions along the error analysis.
In this work, we develop nonstandard and sharp error analysis for the averaging Crank-Nicolson scheme for the subdiffusion, and further propose correction methods to improve the accuracy $O(\tau^\alpha)$ to optimal $O(\tau^2)$ with detailed theoretical analysis.
\par
The rest of the paper is organized as follows.
In Section \ref{sec.prel}, some facts on the approximation of fractional calculus by CQ methods and space semidiscrete schemes for the subdiffusion are provided.
In Section \ref{sec.bdf}, the error analysis for the low accuracy of the averaging Crank-Nicolson scheme is carried out rigorously.
Moreover, correction methods are developed in Section \ref{sec.corr} with detailed error analysis showing that the modified scheme is optimal. 
Several numerical experiments are implemented to verify the theoretical results in Section \ref{sec.expr}.
Finally, some comments are given in Section \ref{sec.conc}.
\par
Throughout the work, by $A \lesssim B$ we mean there exists a positive constant $C$ which may be different at different occurrence such that $A \leq C B$.

\section{Preliminaries}\label{sec.prel}
\subsection{Discrete convolution in CQ}
Denote by $D_t^\alpha$ the Riemann-Liouville fractional differential operator of order $\alpha \in (0,1)$, which is defined by
\begin{equation*}\begin{split}
D_t^\alpha \phi(t)=
\frac{1}{\Gamma(1-\alpha)}
\frac{\mathrm{d}}{\mathrm{d}t}
\int_0^t
\frac{\phi(s)}{(t-s)^\alpha}
\mathrm{d}s.
\end{split}\end{equation*}
The operator $D_t^\alpha$ is closely related to $\partial_t^\alpha$ by the following relation (see \cite{li2015numerical})
\begin{equation*}\begin{split}
\partial_t^\alpha \phi(t)=D_t^\alpha (\phi(t)-\phi(0)).
\end{split}\end{equation*}
\par
On a uniform mesh $0=t_0<t_1<\cdots<t_{N-1}<t_N=T$ with $t_n=n\tau, \tau=T/N$ for some $N>0$, the CQ approximates $D_t^\alpha \phi(t)$ at point $t_{n}$ with the discrete convolution
\begin{equation*}\begin{split}
\mathcal{D}_\tau^\alpha \phi^n=
\tau^{-\alpha}
\sum_{k=0}^n \widetilde{\omega}_{n-k}\phi^k,
\end{split}\end{equation*}
where $\phi^k=\phi(t_k)$ and weights $\widetilde{\omega}_k$ are generated by the function $\widetilde{\omega}(\zeta)$ such that $\widetilde{\omega}(\zeta)=\sum_{k=0}^\infty \widetilde{\omega}_k \zeta^k$.
The CQ theory then asserts (see Theorem 2.5 in \cite{lubich1986discretized}) that $\mathcal{D}_\tau^\alpha \phi^n$ is convergent of order $p$ (to $D_t^\alpha \phi(t_{n})$) if and only if
\begin{equation}\begin{split}\label{pre.0}
\tau^{-\alpha}
\widetilde{\omega}(e^{-\tau})=1+O(\tau^p)\quad \text{and}\quad
\widetilde{\omega}_n=O(n^{-\alpha-1}).
\end{split}\end{equation}
It is clear that the fractional BDF-2 or the generalized 2nd-order Newton–Gregory formula with the following generating functions
\begin{equation}\label{pre.0.1}
\widetilde{\omega}(\zeta)=
\begin{cases}
  \displaystyle \bigg(\frac{3}{2}-2\zeta+\frac{1}{2}\zeta^2 \bigg)^\alpha, & \mbox{for the fractional BDF-2,}\\
  \displaystyle(1-\zeta)^\alpha \bigg[1+\frac{\alpha}{2}(1-\zeta)\bigg], & \mbox{for the generalized 2nd-order Newton–Gregory formula},
\end{cases}
\end{equation}
fulfill the requirements (\ref{pre.0}) with order $p=2$.
By introducing the operator
\begin{equation*}\begin{split}\label{pre.1}
D_\tau^\alpha \phi^{n-\frac{1}{2}}:=
\tau^{-\alpha}
\sum_{k=0}^n \omega_{n-k}\phi^k
\end{split}\end{equation*}
as the approximation to $D_t^\alpha \phi(t)$ at $t=t_{n-\frac{1}{2}}$ and using the averaging technique (\ref{intr.2}), one readily gets
\begin{equation*}\begin{split}
\omega_k=\frac{1}{2}\widetilde{\omega}_k
+\frac{1}{2}\widetilde{\omega}_{k-1}, \quad k\geq 0,\quad \text{with the assumption } \widetilde{\omega}_k=0~ \text{ if ~ $k<0$},
\end{split}\end{equation*}
which leads to the fact that 
\begin{equation}\begin{split}\label{pre.2}
\omega(\zeta)=\frac{1+\zeta}{2}\widetilde{\omega}(\zeta).
\end{split}\end{equation}
\begin{remark}
For generating functions $\widetilde{\omega}(\zeta)$ defined in (\ref{pre.0.1}), one can observe that they are analytic and nonzero on the closed unit disc except $\zeta=1$.
However, the generating function $\omega(\zeta)$ in (\ref{pre.2}) is zero at $\zeta=-1$ (which is on the unit circle) and is therefore problematic if standard error analysis technique such as that in \cite{jin2018analysis} is adopted.
\end{remark}

\subsection{Space semidiscrete scheme and fully discrete scheme}
Since our main interest is on the accuracy of ACN scheme in temporal direction, we simply adopt finite element methods for discretization of spacial variables.
Let $\mathcal{T}_h$ be a shape regular, quasi-uniform triangulation of the domain $\Omega$ where $h$ stands for the mesh size.
Define the space $V_h=\{\chi_h \in H_0^1(\Omega): \chi_h |_e \in \mathcal{P}^1, e \in \mathcal{T}_h\}$
where $\mathcal{P}^1$ denotes the linear polynomial function space.
Introduce the operators  $P_h: L^2(\Omega) \to V_h$ and $R_h: H_0^1(\Omega) \to V_h$ such that
\begin{equation*}\begin{split}
(P_h \phi, \chi_h)=(\phi,\chi_h),\quad \forall \phi \in L^2(\Omega), \forall \chi_h \in V_h,
\\
(\nabla R_h \phi, \nabla \chi_h)=(\nabla \phi, \nabla \chi_h),\quad \forall \phi \in H_0^1(\Omega), \forall \chi_h \in V_h.
\end{split}\end{equation*}
The space semidiscrete scheme is to find $u_h \in V_h$ such that for any $\chi_h \in V_h$, there holds
\begin{equation}\label{pre.3}\begin{split}
(\partial_t^\alpha u_h,\chi_h)+(\nabla u_h,\nabla \chi_h)=(f,\chi_h),
\end{split}\end{equation}
with the initial condition
\begin{equation*}\begin{split}
u_h(0)=v_h:=R_h u_0.
\end{split}\end{equation*}
By further introducing the operator $\Delta_h:V_h \to V_h$ such that
\begin{equation*}\begin{split}
(\Delta_h \phi_h,\chi_h)=-(\nabla \phi_h,\nabla \chi_h),\quad \forall \phi_h,\chi_h \in V_h,
\end{split}\end{equation*}
and letting $w_h:=u_h-v_h$, we rewrite (\ref{pre.3}) as
\begin{equation}\begin{split}\label{pre.4}
D_t^\alpha w_h(t)-\Delta_h w_h(t)=f_h(t)+\Delta_h v_h,\quad t>0,
\end{split}\end{equation}
where $f_h:=P_hf$.
For simplicity, let $f_h(t)=f_h(0)+g_h(t)$ and define $\widetilde{g}_h(\zeta):=\sum_{n=0}^\infty g_h^n \zeta^n$.
It is notable that the operator $\Delta_h$ is sectorial and satisfies the resolvent estimate
\begin{equation}\label{pre.4.1}
\|(z-\Delta_h)^{-1}\| \lesssim |z|^{-1} \quad \text{for any $z\in \Sigma_\sigma$},
\end{equation}
where $\Sigma_\sigma$ denotes the open sector $\{z\in \mathbb{C}: |\arg z|<\sigma, z\neq 0\}$ for some $\sigma \in (\pi/2,\pi)$.
\par
The ACN scheme can be formulated as
\begin{equation}\begin{split}\label{pre.5}
D_\tau^\alpha W_h^{n-\frac{1}{2}}
-\frac{1}{2}\big(
\Delta_h W_h^{n}+\Delta_h W_h^{n-1}\big)
=\frac{1}{2}\big(g_h^{n}+g_h^{n-1}\big)
+f_h(0)
+\Delta_h v_h,
\quad n \geq 1,
\end{split}\end{equation}
or that
\begin{equation*}\begin{split}\label{pre.6}
D_\tau^\alpha W_h^{n-\frac{1}{2}}
-\frac{1}{2}\big(
\Delta_h W_h^{n}+\Delta_h W_h^{n-1}\big)
=\frac{1}{2}\big(f_h^{n}+f_h^{n-1}\big)
+\Delta_h v_h,
\quad n \geq 1.
\end{split}\end{equation*}

\subsection{Contours in the complex plain}
\begin{figure*}
  \includegraphics[width=0.5\textwidth]{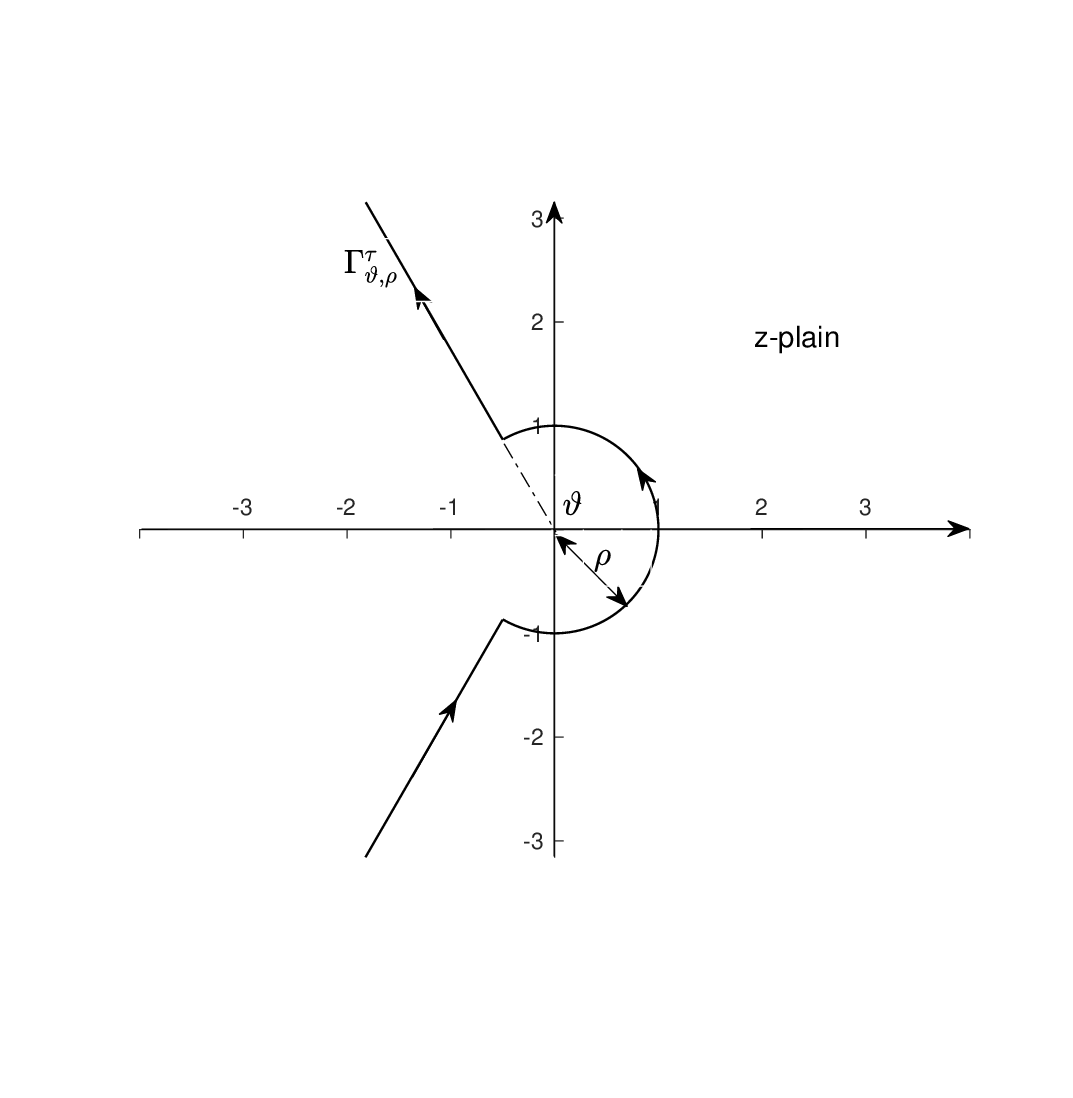}
   \includegraphics[width=0.5\textwidth]{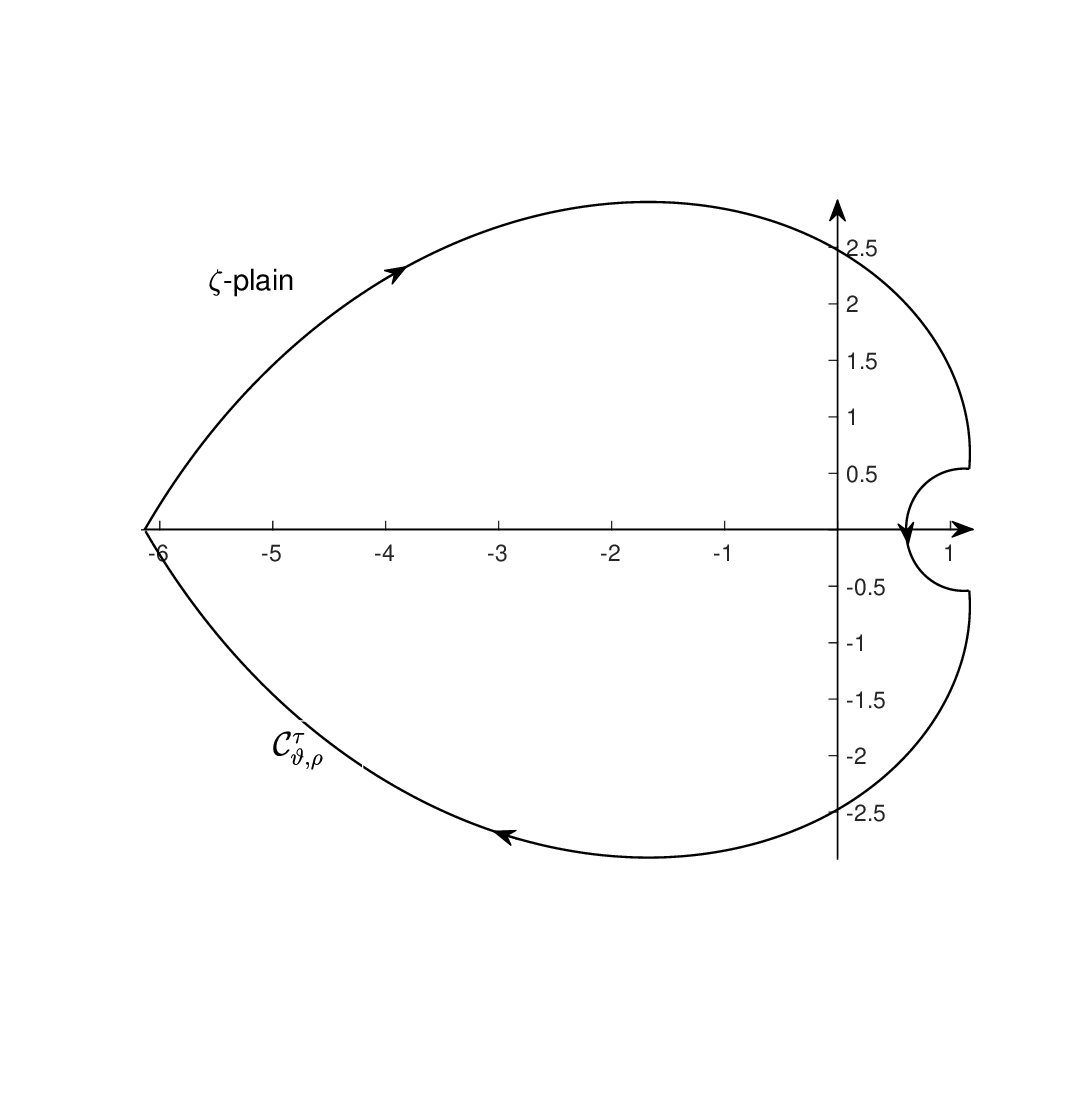}
\caption{Contour $\Gamma_{\vartheta,\rho}^\tau$ in $z$-plain and $\mathcal{C}_{\vartheta,\rho}^\tau$ in $\zeta$-plain under the transformation $\zeta=e^{-z\tau}$.}
\label{fig.1}       
\end{figure*}

\begin{figure*}
  \includegraphics[width=0.5\textwidth]{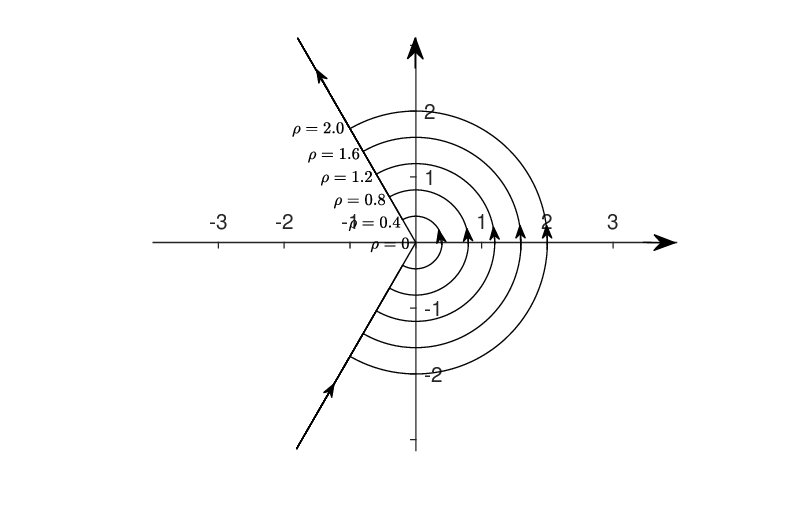}
   \includegraphics[width=0.5\textwidth]{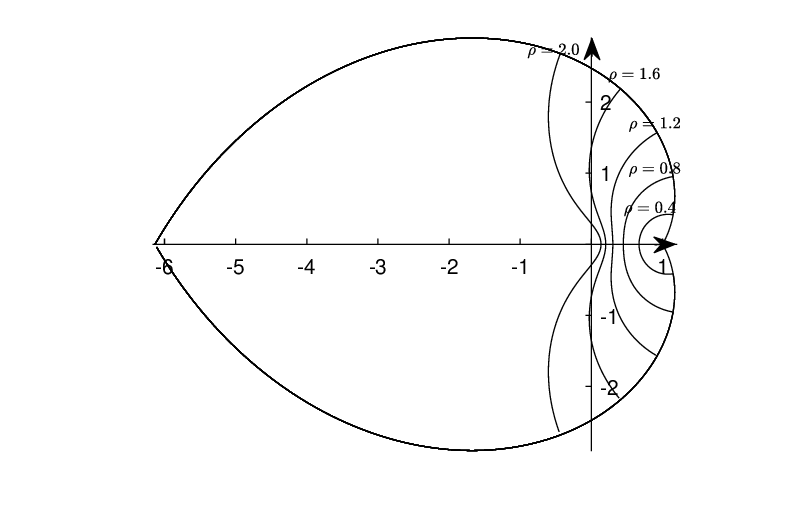}
\caption{Contour $\Gamma_{\vartheta,\rho}^\tau$ in $z$-plain and $\mathcal{C}_{\vartheta,\rho}^\tau$ in $\zeta$-plain under the transformation $\zeta=e^{-z\tau}$ for different $\rho=0,0.4,0.8,1.2,1.6,2.0$, indicating $\mathbb{U}_{\vartheta,\rho_1}^\tau \subset \mathbb{U}_{\vartheta,\rho_2}^\tau$ provided $\rho_1 \geq \rho_2$.}
\label{fig.2}       
\end{figure*}

Let $\Gamma_{\vartheta,\rho}$, $\Gamma_{\vartheta,\rho}^\tau$ be some contours defined by
\begin{equation*}\begin{split}
\displaystyle\Gamma_{\vartheta,\rho}&=\{z\in \mathbb{C}: |z|=\rho, |\arg z|\leq \vartheta\}
\cup
\{z \in \mathbb{C}: z=r e^{\pm {\rm i} \vartheta}, r > \rho\},\quad
\text{for some $\vartheta \in \bigg(\frac{\pi}{2},\pi\bigg)$},
\\
\Gamma_{\vartheta,\rho}^\tau &=
\{z\in \Gamma_{\vartheta,\rho}: |\Im(z)| \leq \pi/\tau\},\quad
\text{for some $\tau \in \bigg(0,\frac{\pi}{\rho\sin\vartheta} \bigg)$},
\end{split}\end{equation*}
oriented with an increasing imaginary part.
The transform $\Gamma_{\vartheta,\rho}^\tau \ni z \mapsto \zeta= e^{-z\tau}$ then convert the contour $\Gamma_{\vartheta,\rho}^\tau$ in $z$-plane into a closed path denoted by $\mathcal{C}_{\vartheta,\rho}^\tau$ in $\zeta$-plane, as illustrated in Fig.\ref{fig.1}.
It is notable that $\mathcal{C}_{\vartheta,\rho}^\tau$ is  oriented clockwise.
Denote by $\mathbb{U}_{\vartheta,\rho}^\tau$ the region enclosed by $\mathcal{C}_{\vartheta,\rho}^\tau$.
\begin{lemma}\label{lem.1}
Given $\tau>0$ sufficiently small and $\vartheta \in (\pi/2,\pi)$.
There hold
\begin{itemize}
\item[{\rm(i)}] The path $\mathcal{C}_{\vartheta,0}^\tau$ (and therefore $\mathbb{U}_{\vartheta,0}^\tau$) is independent of $\tau$.
\item[{\rm(ii)}] If $\rho_1\geq \rho_2 \geq 0$, then $\mathbb{U}_{\vartheta,\rho_1}^\tau \subset \mathbb{U}_{\vartheta,\rho_2}^\tau$.
\end{itemize}
\begin{proof}
By definition we have $\Gamma_{\vartheta,0}^\tau=\{z \in \mathbb{C}: z=r e^{\pm {\rm i} \vartheta}, r \geq 0, r\tau \leq \pi/\sin\vartheta\}$ and
$\xi=e^{-z\tau}=e^{-r\tau(\cos\vartheta\pm \mathrm{i}\sin\vartheta)}$, 
which leads to 
$$\mathcal{C}_{\vartheta,0}^\tau=\{\zeta\in\mathbb{C}:\zeta=e^{-\ell(\cos\vartheta \pm \mathrm{i}\sin\vartheta)}, 0\leq \ell \leq \pi/\sin\vartheta\}$$ and (i) holds immediately.
The inclusion in (ii) can be checked directly in Fig.\ref{fig.2}.
\end{proof}
\end{lemma}
\par
\begin{figure*}
  \includegraphics[width=0.5\textwidth]{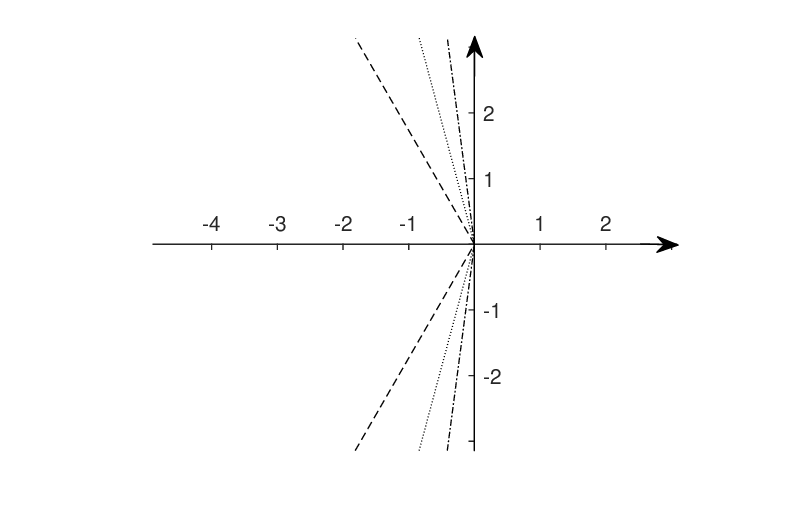}
   \includegraphics[width=0.5\textwidth]{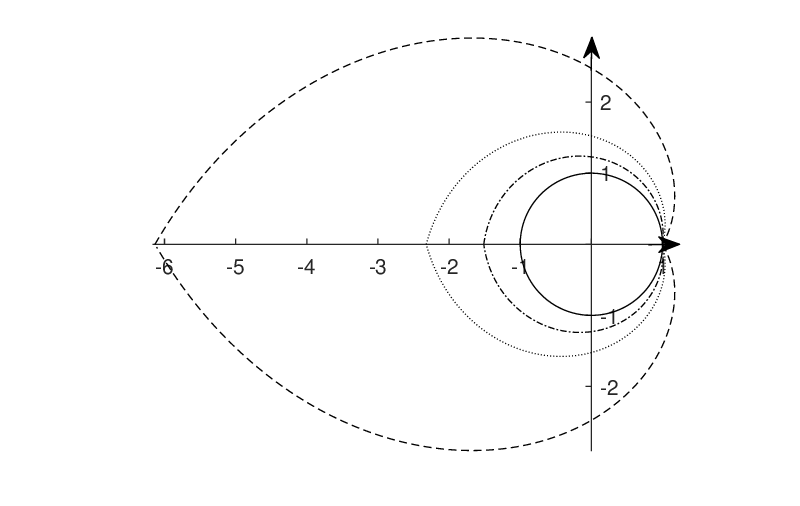}
\caption{Contour $\Gamma_{\vartheta,0}^\tau$ in $z$-plain and $\mathcal{C}_{\vartheta,0}^\tau$ in $\zeta$-plain under the transformation $\zeta=e^{-z\tau}$ for different $\vartheta \in (\pi/2,\pi)$, indicating  $\mathbb{U}_{\vartheta_1} \subset \mathbb{U}_{\vartheta_2}$ provided $\frac{\pi}{2} \leq \vartheta_1 \leq \vartheta_2$.}
\label{fig.3}       
\end{figure*}

The results in (i) permit us to define the path $\mathcal{C}_\vartheta:=\mathcal{C}_{\vartheta,0}^\tau$ (depends only on $\vartheta$) and the region $\mathbb{U}_{\vartheta}:=\mathbb{U}_{\vartheta,0}^\tau$.
By (ii), one gets $\mathbb{U}^\tau_{\vartheta,\rho} \subset \mathbb{U}_\vartheta$ for any $\rho>0$ and $\tau \in (0,\frac{\pi}{\rho\sin\vartheta})$.
More properties of $\mathbb{U}_\vartheta$ are summarized in the next lemma.
\begin{lemma}\label{lem.2}
For the region $\mathbb{U}_{\vartheta}$, there hold
\begin{itemize}
\item[{\rm(i)}] $\mathbb{U}_{\vartheta_1} \subset \mathbb{U}_{\vartheta_2}$ provided $\frac{\pi}{2} \leq \vartheta_1\leq \vartheta_2$,
\item[{\rm(ii)}] $\mathbb{U}_{\frac{\pi}{2}}=\{z: |z|<1\}$,
\item[{\rm(iii)}] For any $z\in \mathbb{C}$ satisfying $|z|>1$, there exits some $\vartheta \in (\pi/2,\pi)$ such that $z \notin \mathbb{U}_\vartheta$,
\item[{\rm(iv)}] For any $\vartheta \in (\pi/2,\pi)$, it holds $-1 \in \mathbb{U}_\vartheta$.
\end{itemize}
\begin{proof}
The inclusion in (i) can be verified directly by Fig. \ref{fig.3}.
For (ii), $\mathbb{U}_{\frac{\pi}{2}}$ is enclosed by the contour $\mathcal{C}_{\frac{\pi}{2}}$ defined by
\begin{equation*}
\{\zeta \in \mathbb{C}: \zeta=e^{\pm {\rm i} \ell}, 0\leq \ell \leq \pi\},
\end{equation*}
which is exactly the unit circle.
(iii) is a natural result of (i) and (ii) since the contour $\mathcal{C}_\vartheta$ narrows continuously to the unit circle as $\vartheta$ tends to $\frac{\pi}{2}$.
For (iv), noting that $-1$ is on the boundary of $\mathbb{U}_{\frac{\pi}{2}}$ by (ii) and using the fact (i), one gets $-1 \in \overline{\mathbb{U}}_\vartheta$ for any $\vartheta \in (\pi/2,\pi)$.
We emphasize that $-1$ can not be on the boundary of $\mathbb{U}_\vartheta$ (the contour $\mathcal{C}_\vartheta$), as $\mathcal{C}_\vartheta$ intersects the real axis only at points $\zeta=1$ (when $z=1$) and $\zeta=-e^{x\tau}$ (when $z=-x+{\rm i}\frac{\pi}{\tau}, x>0$).
\end{proof}
\end{lemma}

\subsection{Solution representation for $w_h(t)$}
Let $\widehat{\phi}(z)$ be the Laplace transform of $\phi(t)$.
Recalling that $\widehat{\partial_t^\alpha \phi}=z^\alpha \widehat{\phi}-z^{\alpha-1}\phi(0)$, we obtain $\widehat{D_t^\alpha w_h(t)}=\widehat{\partial_t^\alpha w_h(t)}=z^\alpha \widehat{w_h}$ thanks to $w_h(0)=0$.
By taking the Laplace transform for (\ref{pre.4}), one gets
\begin{equation*}\label{pre.7}
  (z^\alpha-\Delta_h)\widehat{w_h}=\widehat{g_h}+z^{-1}(\Delta_h v_h+f_h(0)).
\end{equation*}
Resorting to the inverse Laplace transform, we have
\begin{equation}\label{pre.8}
  w_h(t)=\frac{1}{2\pi\mathrm{i}}
  \int_{\Gamma_{\vartheta,\rho}}
  e^{zt}\big(K(z)(\Delta_h v_h+f_h(0))+zK(z)\widehat{g_h}(z)\big)
  \mathrm{d}z,
\end{equation}
where the function $K(z)$ stands for $K(z)=z^{-1}(z^\alpha-\Delta_h)^{-1}$.
\section{Nonstandard error estimates for ACN scheme}\label{sec.bdf}
The standard analysis developed in \cite{jin2018analysis} consists of solution representations of both the time continuous problem (\ref{pre.4}) and the discrete counterpart (\ref{pre.5}) where the latter requires the generating function $\omega(\zeta)$ be analytic and nonzero in the closed unit disc except for the point $\zeta=1$.
Clearly, for the ACN scheme the function $\omega(\zeta)$ in
(\ref{pre.2}) is zero at point $\zeta=-1$ and fails to meet this requirement.
A key observation to surmount this problem is that $\widetilde{\omega}(\zeta)$ for the fractional BDF-2 or the generalized 2nd-order Newton–Gregory formula defined in (\ref{pre.0.1}) satisfies the requirement, and by using the residue theorem and the properties of $\mathbb{U}_\vartheta$ listed in Lemma \ref{lem.2}, one can obtain the expression of the solution immediately.
\begin{lemma}\label{lem.3}
For $\widetilde{\omega}(\zeta)$ defined in (\ref{pre.0.1}), there exists some $\vartheta \in (\pi/2,\pi)$ such that
\end{lemma}
\begin{equation*}\begin{split}
\big\|
\big(
\tau^{-\alpha}
\widetilde{\omega}(\zeta)
-\Delta_h
\big)^{-1}
\big\|
\lesssim
\tau^\alpha 
\big|
\widetilde{\omega}(\zeta)
\big|^{-1},\quad
\forall \zeta \in \mathbb{U}_\vartheta.
\end{split}\end{equation*}
\begin{proof}
In accordance to (\ref{pre.4.1}), we only need to show that $\widetilde{\omega}(\zeta)$ belongs to $\Sigma_\sigma$ for some $\sigma \in (\pi/2,\pi)$.
Indeed, the singular point of $\widetilde{\omega}(\zeta)$ for fractional BDF-2 (except for $\zeta=1$) is $\zeta=3$ and the zero for the generalized 2nd-order Newton–Gregory formula is $\zeta=1+\frac{2}{\alpha}$, both of which satisfy $|\zeta|>1$.
By (iii) in Lemma \ref{lem.2}, there exists some $\vartheta \in (\pi/2,\pi)$ such that $3 \notin \mathbb{U}_\vartheta$ and $1+\frac{2}{\alpha} \notin \mathbb{U}_\vartheta$, meaning that $\widetilde{\omega}(\zeta)$ is analytic and nonzero in $\mathbb{U}_\vartheta$.
Since it is well known that both methods are A$(\frac{\pi}{2})$-stable, i.e., $\widetilde{\omega}(\zeta) \in \Sigma_{\frac{\pi}{2}}$ for any $\zeta \in \mathbb{U}_{\frac{\pi}{2}}=\{z: |z|<1\}$ (see (ii) in Lemma \ref{lem.2}), then by analyticity of $\widetilde{\omega}(\zeta)$ one gets $\widetilde{\omega}(\zeta) \in \Sigma_\sigma$ for some $\sigma \in (\pi/2,\pi)$ so long as $\vartheta$ is close to $\frac{\pi}{2}$ sufficiently.
\end{proof}

\begin{theorem}\label{thm.1}
Given $\tau>0$ and $\alpha \in (0,1)$. 
There exist $\vartheta \in (\pi/2,\pi)$ and $\rho>0$ which are independent of $\tau$ such that the solution of the ACN scheme can be formulated by the following form
\begin{equation}\begin{split}\label{bdf.0}
W_h^n&=
\frac{1}{2\pi {\rm i}}
\int_{\Gamma_{\vartheta,\rho}^\tau}
e^{z t_n} 
\big[
\tau \Phi_i(e^{-z\tau})(\Delta_h v_h+f_h(0))
+\tau\Phi_g(e^{-z\tau})
\widetilde{g}_h(e^{-z\tau})
\big]
\mathrm{d}z
\\&\quad-
(-1)^{n}\big(
\tau^{-\alpha} \widetilde{\omega}(-1)-\Delta_h
\big)^{-1}
\big(
\Delta_h v_h+f_h(0)
\big),
\end{split}\end{equation}
where $\Phi_i(\zeta)$ and $\Phi_g(\zeta)$ are defined by
\begin{equation}\begin{split}\label{bdf.1}
\Phi_i(\zeta)=\frac{2\zeta}{1-\zeta^2}
\big(
\tau^{-\alpha}
\widetilde{\omega}(\zeta)
-\Delta_h
\big)^{-1},\quad
\Phi_g(\zeta)=
\big(
\tau^{-\alpha}
\widetilde{\omega}(\zeta)
-\Delta_h
\big)^{-1}.
\end{split}\end{equation}
\begin{proof}
Multiplying both sides of (\ref{pre.5}) by $\zeta^n$ and summing the index $n$ from $1$ to $\infty$, one gets
\begin{equation}\begin{split}\label{bdf.2}
\bigg(
\tau^{-\alpha}\omega(\zeta)
-\frac{1+\zeta}{2}\Delta_h
\bigg) W_h(\zeta)
=
\frac{1+\zeta}{2}
\widetilde{g}_h(\zeta)
+
\frac{\zeta}{1-\zeta}
\big(f_h(0)+\Delta_h v_h\big),
\end{split}\end{equation}
where we have used the following identity
\begin{equation*}\begin{split}
\sum_{n=1}^\infty 
\zeta^n D_\tau^\alpha W_h^{n-\frac{1}{2}}
=\tau^{-\alpha} \omega(\zeta)W_h(\zeta),
\\
\frac{1}{2}\sum_{n=1}^\infty
\zeta^n
\big(\Delta_h W_h^n+\Delta_h W_h^{n-1}\big)
=\frac{1+\zeta}{2}\Delta_h W_h(\zeta),
\\
\frac{1}{2}\sum_{n=1}^\infty
\zeta^n
(g_h^n+g_h^{n-1})
=\frac{1+\zeta}{2}\widetilde{g}_h(\zeta),
\\
\sum_{n=1}^\infty
\zeta^n
\big(f_h(0)+\Delta_h v_h\big)
=
\frac{\zeta}{1-\zeta}
\big(f_h(0)+\Delta_h v_h\big).
\end{split}\end{equation*}
Clearly, (\ref{bdf.2}) indicates that $\zeta \neq -1$.
By multiplying both sides of (\ref{bdf.2}) with $\frac{2}{1+\zeta}$, we have
\begin{equation*}\begin{split}\label{bdf.2.1}
\big(
\tau^{-\alpha}\widetilde{\omega}(\zeta)
-\Delta_h
\big) W_h(\zeta)
=
\widetilde{g}_h(\zeta)
+
\frac{2\zeta}{1-\zeta^2}
\big(f_h(0)+\Delta_h v_h\big),
\end{split}\end{equation*}
where the operator $\tau^{-\alpha}\widetilde{\omega}(\zeta)-\Delta_h$ is invertible for $\zeta \in \mathbb{U}_\vartheta$ for some $\vartheta \in (\pi/2,\pi)$, due to Lemma \ref{lem.3}.
We therefore obtained the following relation, with $\Phi_i(\zeta)$ and $\Phi_g(\zeta)$ defined in (\ref{bdf.1}), that
\begin{equation}\begin{split}\label{bdf.2.2}
W_h(\zeta)
=
\Phi_i(\zeta)
\big(\Delta_h v_h + f_h(0) \big)
+
\Phi_g(\zeta)
\widetilde{g}_h(\zeta),\quad \forall \zeta \in \mathbb{U}_\vartheta \setminus \{-1\}.
\end{split}\end{equation}
For an appropriately chosen $\rho>0$, the result (ii) in Lemma \ref{lem.1} tells us $\mathbb{U}_{\vartheta,\rho}^\tau \subset \mathbb{U}_{\vartheta,0}^\tau=\mathbb{U}_\vartheta$, implying that $(\ref{bdf.2.2})$ also holds for $\zeta \in \mathbb{U}_{\vartheta,\rho}^\tau \setminus \{-1\}$.
Using the residue theorem, one gets immediately that
\begin{equation}\begin{split}\label{bdf.3}
-\frac{1}{2\pi{\rm i}}\int_{\mathcal{C}_{\vartheta,\rho}^\tau}
\frac{W_h(\zeta)}{\zeta^{n+1}}\mathrm{d}\zeta
={\rm Res}\bigg(\frac{W_h(\zeta)}{\zeta^{n+1}},0\bigg)
+{\rm Res}\bigg(\frac{W_h(\zeta)}{\zeta^{n+1}},-1\bigg),
\end{split}\end{equation}
where the first term of the right side of the equation is exactly $W_h^n$, and the second term can further be formulated as
\begin{equation}\begin{split}\label{bdf.4}
{\rm Res}\bigg(\frac{W_h(\zeta)}{\zeta^{n+1}},-1\bigg) 
&=(-1)^{n-1}
{\rm Res}\big(\Phi_i(\zeta),-1\big)
(\Delta_h v_h+f_h(0))
\\&=
(-1)^{n}\big(
\tau^{-\alpha} \widetilde{\omega}(-1)-\Delta_h
\big)^{-1}
\big(
\Delta_h v_h+f_h(0)
\big).
\end{split}\end{equation}
By setting $\zeta=e^{-z\tau}$, the term of the left side of (\ref{bdf.3}) reads that
\begin{equation*}\begin{split}
-\frac{1}{2\pi{\rm i}}\int_{\mathcal{C}_{\vartheta,\rho}^\tau}
\frac{W_h(\zeta)}{\zeta^{n+1}}\mathrm{d}\zeta=
\frac{1}{2\pi {\rm i}}
\int_{\Gamma_{\vartheta,\rho}^\tau}
e^{z t_n} 
\big[
\tau \Phi_i(e^{-z\tau})(\Delta_h v_h+f_h(0))
+\tau\Phi_g(e^{-z\tau})
\widetilde{g}_h(e^{-z\tau})
\big]
\mathrm{d}z,
\end{split}\end{equation*}
which, combined with (\ref{bdf.4}), completes the proof of the theorem.
\end{proof}
\end{theorem}
The next lemma is from Lemma B.1 in \cite{jin2017correction} and Lemma 2 in \cite{wang2020higher} (Actually, Lemma 2 in \cite{jin2017correction} indicates that (\ref{bdf.4.1}) holds for any $z \in \Gamma_{\vartheta,0}^\tau$.
By introducing a positive $\rho$, one can check readily (\ref{bdf.4.1}) still holds and the proof is omitted here).
\begin{lemma}\label{lem.4}
Given $\tau>0$.
For $\widetilde{\omega}(\zeta)$ defined in (\ref{pre.0.1}),
there exist some $\vartheta' \in (\pi/2,\pi)$ and $\rho>0$ which are independent of $\tau$, such that
\begin{equation}\label{bdf.4.1}
 \tau^\alpha |z|^\alpha
 \lesssim
 |\widetilde{\omega}(e^{-z\tau})|
 \lesssim
 \tau^\alpha |z|^\alpha,\quad \forall \vartheta\in (\pi/2,\vartheta'),\quad \forall z \in \Gamma_{\vartheta,\rho}^\tau.
\end{equation}
\end{lemma}
In combination with Lemma \ref{lem.3} and Lemma \ref{lem.4}, we obtain
\begin{equation}\label{bdf.4.2}
  \big\|\big(\tau^{-\alpha}\widetilde{\omega}(e^{-z\tau})-\Delta_h\big)^{-1} \big\| \lesssim |z|^{-\alpha},\quad \forall z \in \Gamma_{\vartheta,\rho}^\tau,
\end{equation}
where $\vartheta \in (\pi/2,\pi)$ and $\rho>0$ are some constants free of $\tau$.
\begin{lemma}\label{lem.4.1}
Given $\tau>0$.
Let $\widetilde{\omega}(\zeta)$ be defined in (\ref{pre.0.1}).
There exist some $\vartheta' \in (\pi/2,\pi)$ and $\rho>0$ which are independent of $\tau$, such that
\begin{equation}\label{bdf.4.2.1}
|z^\alpha-\tau^{-\alpha}\widetilde{\omega}(e^{-z\tau})|\lesssim
\tau^2|z|^{2+\alpha},\quad \forall \vartheta\in (\pi/2,\vartheta'),\quad
\forall z \in \Gamma_{\vartheta,\rho}^\tau.
\end{equation}
\begin{proof}
For the fractional BDF-2, Lemma B.1 in \cite{jin2017correction} has shown (\ref{bdf.4.2.1}) holds.
The technique can also be applied to the GNG-2 similarly and the details are omitted here.
\end{proof}
\end{lemma}
\begin{lemma}\label{lem.5}
Given $\tau>0$, there exist some $\vartheta \in (\pi/2,\pi)$ and $\rho>0$ which are independent of $\tau$ such that
\begin{itemize}
\item[{\rm(i)}]
$\displaystyle\bigg|
\frac{2e^{-z\tau}}{1-e^{-2z\tau}}-z^{-1}\tau^{-1}
\bigg| \lesssim \tau|z|,
\quad
\bigg|
\frac{e^{-z\tau}}{(1-e^{-z\tau})^2}
-z^{-2}\tau^{-2}
\bigg|\lesssim 1,
\quad \forall z\in \Gamma_{\vartheta,\rho}^\tau$.
\item[{\rm(ii)}]
$\|\tau \Phi_i(e^{-z\tau})-K(z)\|
\lesssim \tau^2|z|^{1-\alpha},
\quad \forall z\in \Gamma_{\vartheta,\rho}^\tau$.
\item[{\rm(iii)}]
$\displaystyle
\bigg\|
\frac{\tau^2 e^{-z\tau}}{(1-e^{-z\tau})^2}
\Phi_g(e^{-z\tau})
-z^{-1}K(z)
\bigg\|
\lesssim \tau^2|z|^{-\alpha},
\quad \forall z\in \Gamma_{\vartheta,\rho}^\tau,
$
\end{itemize}
where $K(z)=z^{-1}(z^\alpha-\Delta_h)^{-1}$ and $\Phi_i(\zeta)$ and $\Phi_g(\zeta)$ are defined in (\ref{bdf.1}).
\begin{proof}
For (i), if $|z|\tau<\epsilon$ where $\epsilon>0$ is given small enough, by using the Taylor expansion of $e^{-z\tau}$, one can check (i) holds readily.
To be specific, denote by $Q_n(z)$ some function which may not be the same at different occurrence, fulfilling $|Q_n(z)|\lesssim |z|^n$ for sufficiently small $|z|$.
There holds
\begin{equation*}\begin{split}\label{bdf.4.2.1.1}
\frac{2z\tau e^{-z\tau}}{1-e^{-2z\tau}}
&=\frac{2z\tau\big(1-z\tau+z^2\tau^2/2+Q_3(z\tau)\big)}
{1-\big(1-2z\tau+2z^2\tau^2+Q_3(z\tau)\big)}
=\frac{1-z\tau+z^2\tau^2/2+Q_3(z\tau)}{1-z\tau-Q_2(z\tau)}
\\
&=
1+\frac{z^2\tau^2/2+Q_3(z\tau)+Q_2(z\tau)}{1-z\tau-Q_2(z\tau)}
=1+Q_2(z\tau),
\end{split}
\end{equation*}
which leads to
\begin{equation*}
\bigg|
\frac{2e^{-z\tau}}{1-e^{-2z\tau}}-z^{-1}\tau^{-1}
\bigg|
=|z|^{-1}\tau^{-1}
\bigg|
\frac{2z\tau e^{-z\tau}}{1-e^{-2z\tau}}-1
\bigg|
=|z|^{-1}\tau^{-1}Q_2(z\tau)
\lesssim \tau|z|.
\end{equation*}
To prove the second inequality of (i), we note that
\begin{equation*}\begin{split}
\frac{z^2\tau^2 e^{-z\tau}}{(1-e^{-z\tau})^2}
&=\frac{1-z\tau+z^2\tau^2/2+Q_3(z\tau)}{(1-z\tau/2-Q_2(z\tau))^2}
=\frac{1-z\tau+z^2\tau^2/2+Q_3(z\tau)}{1-z\tau+Q_2(z\tau)}
\\&=
1+\frac{Q_2(z\tau)}{1-z\tau+Q_2(z\tau)},
\end{split}
\end{equation*}
yielding
\begin{equation*}
\bigg|
\frac{e^{-z\tau}}{(1-e^{-z\tau})^2}
-z^{-2}\tau^{-2}
\bigg|
=|z|^{-2}\tau^{-2}
\bigg|
\frac{z^2\tau^2 e^{-z\tau}}{(1-e^{-z\tau})^2}-1
\bigg|
\lesssim 1.
\end{equation*}
\par
If $|z|\tau \geq \epsilon$, let $z=re^{\pm {\rm i}\vartheta}$ with $r \in [\frac{\epsilon}{\tau},\frac{\pi}{\tau\sin\vartheta}]$ for some $\vartheta \in (\pi/2,\pi)$ to be determined next.
We only need to prove that the left terms of the inequalities in (i) can be bounded by some positive constants (independent of $\tau$), since any positive constant satisfies $c \leq \frac{c}{\epsilon}\epsilon \leq \frac{c}{\epsilon}|z|\tau$ by the assumption.
\par
Let $\vartheta',\vartheta'' \in (\pi/2,\pi)$ be given and sufficiently close to $\pi/2$ with $\vartheta'<\vartheta''$.
For any $\vartheta \in (\vartheta',\vartheta'')$, using the following estimates,
\begin{equation*}\begin{split}
|e^{-z\tau}|&\leq e^{|z|\tau}\leq e^{\frac{\pi}{\sin\vartheta}}
\leq 
e^{\frac{\pi}{\sin\vartheta''}},\quad
|z^{-1}\tau^{-1}| \leq \frac{1}{\epsilon},
\\
\big|1-e^{-z\tau}\big|^2
&=(1-e^{-\tau r \cos\vartheta})^2
+2e^{-\tau r \cos 2\vartheta}
[1-\cos(\tau r \sin \vartheta)]
\\
&\geq
(e^{\tau r \sin(\vartheta-\pi/2)}-1)^2
\geq \tau^2 r^2 \sin^2(\vartheta-\pi/2)
\\
&\geq \epsilon^2\sin^2(\vartheta'-\pi/2),
\\
\big|1-e^{-2z\tau}\big|^2&\geq
4\tau^2 r^2\sin^2(\vartheta-\pi/2)
\geq 4\epsilon^2\sin^2(\vartheta'-\pi/2),
\end{split}
\end{equation*}
one can get
\begin{equation*}
\bigg|
\frac{2e^{-z\tau}}{1-e^{-2z\tau}}-z^{-1}\tau^{-1}
\bigg|
\lesssim
\frac{|e^{-z\tau}|}{|1-e^{-2z\tau}|}
+\frac{1}{|z\tau|}
\lesssim
\frac{e^{\frac{\pi}{\sin\vartheta''}}}{2\epsilon\sin(\vartheta'-\pi/2)}
+\frac{1}{\epsilon},
\end{equation*}
and
\begin{equation*}
\bigg|
\frac{e^{-z\tau}}{(1-e^{-z\tau})^2}
-z^{-2}\tau^{-2}
\bigg|
\lesssim
\frac{|e^{-z\tau}|}{|1-e^{-z\tau}|^2}
+\frac{1}{|z\tau|^2}
\lesssim
\frac{e^{\frac{\pi}{\sin\vartheta''}}}{\epsilon^2\sin^2(\vartheta'-\pi/2)}
+\frac{1}{\epsilon^2}.
\end{equation*}
\par
For (ii), we have, by the definition of $\Phi_i$ and $K(z)$ that
\begin{equation}\begin{split}\label{bdf.4.2.2}
\tau\Phi_i(e^{-z\tau})-K(z)
&=\tau\bigg(
\frac{2e^{-z\tau}}{1-e^{-2z\tau}}
-z^{-1}\tau^{-1}
\bigg)
\big(
\tau^{-\alpha}
\widetilde{\omega}(e^{-z\tau})
-\Delta_h
\big)^{-1}
\\
&\quad+
z^{-1}
\big[
\big(
\tau^{-\alpha}\widetilde{\omega}(e^{-z\tau})-\Delta_h
\big)^{-1}
-(z^\alpha-\Delta_h)^{-1}
\big].
\end{split}
\end{equation}
In accordance to the following identity
\begin{equation}\begin{split}\label{bdf.4.2.3}
&\quad\big(
\tau^{-\alpha}\widetilde{\omega}(e^{-z\tau})-\Delta_h
\big)^{-1}
-(z^\alpha-\Delta_h)^{-1}
\\&=
\big(
\tau^{-\alpha}\widetilde{\omega}(e^{-z\tau})-\Delta_h
\big)^{-1}
\big(z^\alpha-\tau^{-\alpha}\widetilde{\omega}(e^{-z\tau})\big)
(z^\alpha-\Delta_h)^{-1},
\end{split}
\end{equation}
and the estimates in (i), (\ref{bdf.4.2}), (\ref{bdf.4.2.1}) and (\ref{pre.4.1}), by (\ref{bdf.4.2.2}) one immediately gets
\begin{equation*}\label{bdf.4.2.4}
\|\tau\Phi_i(e^{-z\tau})-K(z)\|\lesssim \tau^2|z|^{1-\alpha}.
\end{equation*}
\par
For (iii), with $\Phi_g$ defined in (\ref{bdf.1}), we have
\begin{equation*}\begin{split}\label{bdf.4.2.5}
\frac{\tau^2 e^{-z\tau}}{(1-e^{-z\tau})^2}
\Phi_g(e^{-z\tau})
-z^{-1}K(z)
&=
\tau^2\bigg[
\frac{e^{-z\tau}}{(1-e^{-z\tau})^2}
-z^{-2}\tau^{-2}
\bigg]
\big(
\tau^{-\alpha}\widetilde{\omega}(e^{-z\tau})
-\Delta_h
\big)^{-1}
\\&\quad
+z^{-2}
\big[
\big(
\tau^{-\alpha}\widetilde{\omega}(e^{-z\tau})
-\Delta_h
\big)^{-1}
-(z^\alpha-\Delta_h)^{-1}
\big],
\end{split}
\end{equation*}
which, in combination with the identity (\ref{bdf.4.2.3}) and (i), (\ref{bdf.4.2}), (\ref{bdf.4.2.1}) and (\ref{pre.4.1}), yields  
\begin{equation*}
\bigg\|
\frac{\tau^2 e^{-z\tau}}{(1-e^{-z\tau})^2}
\Phi_g(e^{-z\tau})
-z^{-1}K(z)
\bigg\|
\lesssim \tau^2|z|^{-\alpha}.
\end{equation*}
\end{proof}
\end{lemma}
\par
Since the source term $f(x,t)$ can be expanded as $f(x,t)=q_0(x)+tq_1(x)+t*q(x,t)$, the sharp error analysis for the ACN scheme consists of the following several theorems.
\begin{theorem}\label{thm.2}
Assume $f(x,t)\equiv 0$.
Let $W_h^n$ be the solution of (\ref{pre.5}) and $U_h^n:=W_h^n+u_h(0)$ be the approximation to $u_h(t_n)$ for $n\geq 1$.
For sufficiently small $\tau>0$, there holds
\begin{equation*}\label{bdf.5}
  \|u_h(t_n)-U_h^n\|_{L^2(\Omega)}\lesssim (t_n^{\alpha-2}\tau^2+\tau^\alpha) \|\Delta u_0\|_{L^2(\Omega)}.
\end{equation*}
\begin{proof}
By definition we have $u_h(t_n)-U_h^n=w_h(t_n)-W_h^n$.
Using the solution representations (\ref{pre.8}) and (\ref{bdf.0}), there holds
\begin{equation*}\begin{split}\label{bdf.4.3}
 w_h(t_n)-W_h^n=I_1^n+I_2^n+I_3^n,
\end{split}\end{equation*}
where $I_i^n (i=1,2,3)$ stand for
\begin{equation*}\begin{split}\label{bdf.4.4}
I_1^n&=\frac{1}{2\pi{\rm i}}
\int_{\Gamma_{\vartheta,\rho}^\tau}e^{zt_n}
\big[K(z)-\tau\Phi_i(e^{-z\tau})\big]
\Delta_h v_h \mathrm{d}z,
\\
I_2^n&=\frac{1}{2\pi{\rm i}}
\int_{\Gamma_{\vartheta,\rho} \setminus \Gamma_{\vartheta,\rho}^\tau}e^{zt_n}
K(z)
\Delta_h v_h \mathrm{d}z,
\\
I_3^n&=(-1)^{n}\big(
\tau^{-\alpha} \widetilde{\omega}(-1)-\Delta_h
\big)^{-1}
\Delta_h v_h.
\end{split}\end{equation*}
For the term $I_1^n$, by appealing to the estimate (ii) in Lemma \ref{lem.5} and the symmetry of contour $\Gamma_{\vartheta,\rho}^\tau$, we have
\begin{equation}\begin{split}\label{bdf.4.5}
\|I_1^n\|_{L^2(\Omega)}
&\lesssim \tau^2\|\Delta_h v_h\|_{L^2(\Omega)}
\int_{\Gamma_{\vartheta,\rho}^\tau}|e^{z t_n}||z|^{1-\alpha}
|\mathrm{d}z|
\\ 
&\lesssim\tau^2\|\Delta_h v_h\|_{L^2(\Omega)}
\bigg(
\int_{\rho}^{\frac{\pi}{\tau\sin\vartheta}}
e^{r t_n \cos\vartheta} r^{1-\alpha}
\mathrm{d}r
+\rho^{2-\alpha}
\int_{0}^{\vartheta}
e^{\rho t_n\cos\theta}
\mathrm{d}\theta
\bigg).
\end{split}
\end{equation}
Taking $s=r t_n$, the first integration in (\ref{bdf.4.5}) reads that
\begin{equation}\label{bdf.4.6}
\int_{\rho}^{\frac{\pi}{\tau\sin\vartheta}}
e^{rt_n\cos\vartheta}r^{1-\alpha}
\mathrm{d}r
\leq
t_n^{\alpha-2}
\int_{0}^{+\infty}e^{s\cos\vartheta}s^{1-\alpha}\mathrm{d}s
\lesssim t_n^{\alpha-2},
\end{equation}
where the last inequality holds for that $\vartheta \in (\pi/2,\pi)$ and $e^{s\cos\vartheta}$ decays faster than $s^{1-\alpha}$ at infinity.
By choosing $\rho$ such that $\rho t_n \leq 1$, the second integration in (\ref{bdf.4.5}) can further be formulated as
\begin{equation*}\label{bdf.4.7}
\rho^{2-\alpha}
\int_{0}^{\vartheta}
e^{\rho t_n \cos\theta}
\mathrm{d}\theta
\leq
t_n^{\alpha-2}\int_{0}^{\vartheta}e\mathrm{d}\theta
\lesssim t_n^{\alpha-2},
\end{equation*}
which, in combination with (\ref{bdf.4.6}) and (\ref{bdf.4.5}), leads to
\begin{equation*}\label{bdf.4.8}
\|I_1^n\|_{L^2(\Omega)}\lesssim t_n^{\alpha-2}\tau^2\|\Delta_h v_h\|_{L^2(\Omega)}.
\end{equation*}
For the term $I_2^n$, since $\|K(z)\|=\|z^{-1}(z^\alpha-\Delta_h)^{-1}\| \lesssim |z|^{-1-\alpha}$, we obtain
\begin{equation}\begin{split}\label{bdf.4.9}
\|I_2^n\|_{L^2(\Omega)}
&\lesssim
\|\Delta_h v_h\|_{L^2(\Omega)}
\int_{\frac{\pi}{\tau\sin\vartheta}}^{+\infty}
e^{r t_n \cos\vartheta}
r^{-1-\alpha}\mathrm{d}r
\lesssim
\tau^2 \|\Delta_h v_h\|_{L^2(\Omega)}
\int_{\frac{\pi}{\tau\sin\vartheta}}^{+\infty}
e^{r t_n \cos\vartheta}
r^{1-\alpha}\mathrm{d}r
\\&
\lesssim
t_n^{\alpha-2}\tau^2\|\Delta_h v_h\|_{L^2(\Omega)},
\end{split}
\end{equation}
where the last inequality is worked out by setting $s=rt_n$ and the fast decay of $e^{s\cos\vartheta}$ as $s$ tends to infinity.
\\
For the term $I_3^n$, direct estimate shows that
\begin{equation}\label{bdf.4.10}
\|I_3^n\|_{L^2(\Omega)}\lesssim \tau^\alpha \|\Delta_h v_h\|_{L^2(\Omega)}.
\end{equation}
Combining (\ref{bdf.4.5}), (\ref{bdf.4.9}) with (\ref{bdf.4.10}) and using the identity $\Delta_h R_h=P_h \Delta$, we complete the proof of the theorem.
\end{proof}
\end{theorem}

\begin{theorem}\label{thm.3}
Assume $u_0(x)\equiv 0$ and $f(x,t)=q_0(x)+tq_1(x)$ with $q_0,q_1 \in L^2(\Omega)$.
Let $W_h^n$ be the solution of (\ref{pre.5}) and $U_h^n:=W_h^n+u_h(0)$ be the approximation to $u_h(t_n)$ for $n\geq 1$.
For sufficiently small $\tau>0$, there holds
\begin{equation*}\label{bdf.6}
  \|u_h(t_n)-U_h^n\|_{L^2(\Omega)}\lesssim 
  \tau^2
  \big(t_n^{\alpha-2}\|q_0(x)\|_{L^2(\Omega)}+t_n^{\alpha-1}\|q_1(x)\|_{L^2(\Omega)}
  \big)
  +\tau^\alpha \|q_0(x)\|_{L^2(\Omega)}.
\end{equation*}
\begin{proof}
For $f(x,t)=q_0(x)+tq_1(x)$, we have $f_h(0)=q_{0h}(x):=P_h q_0(x)$, $g_h(x,t)=tq_{1h}(x):=t P_h q_1(x)$ and
\begin{equation*}\label{bdf.6.1}
 \widetilde{g}_h(\zeta)=q_{1h}(x)\sum_{i=0}^{\infty}t_i \zeta^i
 =\frac{\tau\zeta}{(1-\zeta)^2}q_{1h}(x),
 \quad
 \widehat{g_h}(z)=z^{-2}q_{1h}(x).
\end{equation*}
It is notable that if $q_1(x) \equiv 0$, by almost the same estimate as Theorem \ref{thm.2}, one can obtain
\begin{equation*}
\|u_h(t_n)-U_h^n\|_{L^2(\Omega)}
\lesssim
(t_n^{\alpha-2}\tau^2+\tau^\alpha)\|q_0(x)\|_{L^2(\Omega)}.
\end{equation*}
We therefore in the next analysis mainly focus on $q_0(x) \equiv 0$.
By solution representations (\ref{pre.8}) and (\ref{bdf.0}), we have
\begin{equation}\begin{split}\label{bdf.6.2}
 w_h(t_n)-W_h^n=J_1^n+J_2^n,
\end{split}\end{equation}
where $J_i^n (i=1,2)$ stand for
\begin{equation*}\begin{split}
J_1^n &=
\frac{1}{2\pi {\rm i}}
\int_{\Gamma_{\vartheta,\rho}^\tau}
e^{z t_n}
\big(zK(z)\widehat{g_h}(z)
-\tau\Phi_g(e^{-z\tau})
\widetilde{g}_h(e^{-z\tau})\big)
\mathrm{d}z,
\\
J_2^n &=
\frac{1}{2\pi{\rm i}}
\int_{\Gamma_{\vartheta,\rho}\setminus \Gamma_{\vartheta,\rho}^\tau}
e^{zt_n}
zK(z)\widehat{g_h}(z)
\mathrm{d}z.
\end{split}
\end{equation*}
For the term $J_1^n$, by resorting to the estimate (iii) in Lemma \ref{lem.5} and the symmetry of contour $\Gamma_{\vartheta,\rho}^\tau$, we have
\begin{equation}\begin{split}\label{bdf.6.3}
 \|J_1^n\|_{L^2(\Omega)}
 &\lesssim
 \tau^2\|q_{1h}(x)\|_{L^2(\Omega)}
\bigg(
\int_{\rho}^{\frac{\pi}{\tau\sin\vartheta}}
e^{r t_n \cos\vartheta} r^{-\alpha}
\mathrm{d}r
+\rho^{1-\alpha}
\int_{0}^{\vartheta}
e^{\rho t_n\cos\theta}
\mathrm{d}\theta
\bigg)
\\&
\lesssim
t_n^{\alpha-1}\tau^2\|q_{1h}(x)\|_{L^2(\Omega)},
\end{split}
\end{equation}
where the replacement $s=rt_n$ and assumption $\rho t_n \leq 1$ are adopt in deriving the last inequality. 
\\
For the term $J_2^n$, since $\|zK(z)\widehat{g_h}(z)\|_{L^2(\Omega)}
=\|z^{-1}K(z)q_{1h}(x)\|_{L^2(\Omega)}
\lesssim
|z|^{-2-\alpha}
\|q_{1h}(x)\|_{L^2(\Omega)}$, one gets
\begin{equation}\begin{split}\label{bdf.6.4}
\|J_2^n\|_{L^2(\Omega)} 
&\lesssim
\|q_{1h}(x)\|_{L^2(\Omega)}
\int_{\frac{\pi}{\tau\sin\vartheta}}^{+\infty}
e^{r t_n \cos\vartheta}
r^{-2-\alpha}\mathrm{d}r
\lesssim
\tau^2 \|\Delta_h v_h\|_{L^2(\Omega)}
\int_{\frac{\pi}{\tau\sin\vartheta}}^{+\infty}
e^{r t_n \cos\vartheta}
r^{-\alpha}\mathrm{d}r
\\&
\lesssim
t_n^{\alpha-1}\tau^2\|q_{1h}(x)\|_{L^2(\Omega)}.
\end{split}
\end{equation}
Combining (\ref{bdf.6.2}), (\ref{bdf.6.3}) with (\ref{bdf.6.4}) and using the fact that $\|q_{1h}(x)\|_{L^2(\Omega)} \leq \|q_1(x)\|_{L^2(\Omega)}$, we get
\begin{equation*}
\|w_h(t_n)-W_h^n\|_{L^2(\Omega)}
\lesssim
t_n^{\alpha-1}
\tau^2\|q_{1h}(x)\|_{L^2(\Omega)}.
\end{equation*}
The proof of the theorem is completed.
\end{proof}
\end{theorem}
\begin{remark}
Theorem \ref{thm.3} indicates that if the source term $f(x,t)$ is smooth and is zero at time $t=0$, such as $f=tq_1(x)$, then the error at $t_n$ is of $O(\tau^2)$.
See Example 2 in Section \ref{sec.expr} for verification. 
\end{remark}

\begin{theorem}\label{thm.4}
Assume $u_0(x)\equiv 0$ and $f(x,t)=t*q(x,t)$.
Let $W_h^n$ be the solution of (\ref{pre.5}) and $U_h^n:=W_h^n+u_h(0)$ be the approximation to $u_h(t_n)$ for $n\geq 1$.
For sufficiently small $\tau>0$, there holds
\begin{equation}\label{bdf.7}
  \|u_h(t_n)-U_h^n\|_{L^2(\Omega)}
  \lesssim 
  \tau^2\int_{0}^{t_n}(t_n-s)^{\alpha-1}\|q(s)\|_{L^2(\Omega)}\mathrm{d}s.
\end{equation}
\begin{proof}
In this case, the discrete solution representation in (\ref{bdf.0}) consists only the integral part.
The estimate (\ref{bdf.7}) can be carried out using the same analysis as Lemma 3.12 in \cite{jin2018analysis} which is omitted here.
\end{proof}
\end{theorem}
\begin{remark}
Combining Theorem \ref{thm.2}-Theorem \ref{thm.4}, we actually have shown the following result that,
by assuming $u_0(x) \in D(\Delta)$, $f(x,t) \in W^{1,\infty}(0,T;L^2(\Omega))$, $\int_{0}^{t}(t-s)^{\alpha-1}\|f''(s)\|_{L^2(\Omega)}\mathrm{d}s<\infty$,
and letting $W_h^n$ be the solution of (\ref{pre.5}) and $U_h^n:=W_h^n+u_h(0)$ be the approximation to $u_h(t_n)$ for $n\geq 1$,
then for sufficiently small $\tau>0$, there holds
\begin{equation*}\begin{split}\label{bdf.8}
  \|u_h(t_n)-U_h^n\|_{L^2(\Omega)}
  &\lesssim 
  \tau^2
  \bigg(
  t_n^{\alpha-2}(\|f(0)\|_{L^2(\Omega)}+\|\Delta u_0\|_{L^2(\Omega)})
  +t_n^{\alpha-1}\|f'(0)\|_{L^2(\Omega)}
  \\&\quad
  +\int_{0}^{t_n}(t_n-s)^{\alpha-1}\|f''(s)\|_{L^2(\Omega)}\mathrm{d}s
  \bigg)
  +\tau^\alpha(\|\Delta u_0\|_{L^2(\Omega)}+\|f(0)\|_{L^2(\Omega)}).
\end{split}
\end{equation*}
\end{remark}
\section{Correction method for ACN scheme}\label{sec.corr}
In this section we derive the modified averaging Crank-Nicolson (MACN) scheme by adding corrections to obtain the optimal accuracy.
In general, at least two time steps should be modified (in which case at least two parameters are involved) since the singularity of $\Phi_i(\zeta)$ at $\zeta=-1$ must be removed and some estimate, such as (\ref{cor.3}) in the following analysis, should be satisfied.
Define $a^{(1)}$ and $a^{(2)}$ as the correction parameters, with which the modified ACN scheme can be formulated as
\begin{equation}\begin{split}\label{cor.1}
D_\tau^\alpha W_h^{\frac{1}{2}}
-\frac{1}{2}
\Delta_h W_h^{1}
=\frac{1}{2}\big(f_h^{1}+f_h(0)\big)
+\Delta_h v_h
+a^{(1)}(f_h(0)+\Delta_h v_h),
\quad & n=1,
\\
D_\tau^\alpha W_h^{\frac{3}{2}}
-\frac{1}{2}\big(
\Delta_h W_h^{2}+\Delta_h W_h^{1}\big)
=\frac{1}{2}\big(f_h^{2}+f_h^{1}\big)
+\Delta_h v_h
+a^{(2)}(f_h(0)+\Delta_h v_h), \quad & n=2,
\\
D_\tau^\alpha W_h^{n-\frac{1}{2}}
-\frac{1}{2}\big(
\Delta_h W_h^{n}+\Delta_h W_h^{n-1}\big)
=\frac{1}{2}\big(f_h^{n}+f_h^{n-1}\big)
+\Delta_h v_h,
\quad & n \geq 3.
\end{split}\end{equation}
\begin{theorem}
Given $\tau>0$ and $\alpha \in (0,1)$. 
There exist $\vartheta \in (\pi/2,\pi)$ and $\rho>0$ which are independent of $\tau$, such that the solution of the modified ACN scheme (for some carefully chosen $a^{(1)}$ and $a^{(2)}$) takes the form
\begin{equation*}\begin{split}
W_h^n&=
\frac{1}{2\pi {\rm i}}
\int_{\Gamma_{\vartheta,\rho}^\tau}
e^{z t_n} 
\big[
\tau \Phi_{mi}(e^{-z\tau})(\Delta_h v_h+f_h(0))
+\tau\Phi_g(e^{-z\tau})
\widetilde{g}_h(e^{-z\tau})
\big]
\mathrm{d}z
\end{split}\end{equation*}
where $\Phi_g(\zeta)$ is defined in (\ref{bdf.1}) and $\Phi_{mi}(\zeta)$ is defined by
\begin{equation*}\begin{split}
\Phi_{mi}(\zeta)=
\frac{2\zeta}{1+\zeta}
\bigg(
\frac{1}{1-\zeta}+a^{(1)}
+a^{(2)}\zeta
\bigg)
\big(\tau^{-\alpha}\widetilde{\omega}(\zeta)-\Delta_h\big)^{-1}.
\end{split}\end{equation*}
\begin{proof}
Multiplying both sides of (\ref{cor.1}) by $\zeta^n$ and summing the index $n$ from $1$ to $\infty$, one gets
\begin{equation*}\begin{split}
\big(\tau^{-\alpha}\widetilde{\omega}(\zeta)-\Delta_h\big)W_h(\zeta)&=
\widetilde{g}_h(\zeta)
+
\frac{2\zeta}{1+\zeta}
\bigg(
\frac{1}{1-\zeta}+a^{(1)}
+a^{(2)}\zeta
\bigg)
\big(
f_h(0)+\Delta_h v_h
\big),
\end{split}\end{equation*}
which leads to
\begin{equation}\begin{split}\label{cor.1.1}
  W_h(\zeta)&=\frac{2\zeta}{1+\zeta}
\bigg(
\frac{1}{1-\zeta}+a^{(1)}
+a^{(2)}\zeta
\bigg)
\big(\tau^{-\alpha}\widetilde{\omega}(\zeta)-\Delta_h\big)^{-1}
\big(
f_h(0)+\Delta_h v_h
\big)
\\
&\quad+
\big(\tau^{-\alpha}\widetilde{\omega}(\zeta)-\Delta_h\big)^{-1}
\widetilde{g}_h(\zeta)
\\
&=
\Phi_{mi}(\zeta)\big(
f_h(0)+\Delta_h v_h
\big)
+\Phi_g(\zeta)\widetilde{g}_h(\zeta).
\end{split}\end{equation}
For carefully chosen $a^{(1)}$ and $a^{(2)}$ such that the term $\frac{2\zeta}{1+\zeta}
\big(
\frac{1}{1-\zeta}+a^{(1)}
+a^{(2)}\zeta
\big)$
is analytic at $\zeta=-1$, i.e., by requiring that
\begin{equation}\begin{split}\label{cor.1.1.1}
\frac{1}{2}+a^{(1)}-a^{(2)}=0,
\end{split}\end{equation}
we conclude (\ref{cor.1.1}) holds for any $\zeta \in \mathbb{U}^\tau_{\vartheta,\rho}$ for some $\vartheta \in (\pi/2,\pi)$ and $\rho>0$.
Using the residue theorem, we obtain that
\begin{equation}\begin{split}\label{cor.1.2}
-\frac{1}{2\pi{\rm i}}\int_{\mathcal{C}_{\vartheta,\rho}^\tau}
\frac{W_h(\zeta)}{\zeta^{n+1}}\mathrm{d}\zeta
={\rm Res}\bigg(\frac{W_h(\zeta)}{\zeta^{n+1}},0\bigg)=W_h^n.
\end{split}\end{equation}
Then, by setting $\zeta=e^{-z\tau}$ in (\ref{cor.1.2}) we completes the proof of the theorem.
\end{proof}
\end{theorem}
\par
Considering (\ref{cor.1.1.1}), $\Phi_{mi}$ can be reformulated as
\begin{equation*}\begin{split}
\Phi_{mi}(\zeta)=
\frac{2\zeta}{1-\zeta}
\bigg[
1+a^{(1)}
-\bigg(a^{(1)}+\frac{1}{2}\bigg)\zeta
\bigg]
\big(\tau^{-\alpha}\widetilde{\omega}(\zeta)-\Delta_h\big)^{-1}.
\end{split}\end{equation*}
To determine $a^{(1)}$, one needs to justify the estimate (counterpart of (ii) in Lemma \ref{lem.5})
\begin{equation}\begin{split}\label{cor.3}
\big\|\tau \Phi_{mi}(e^{-z\tau})-z^{-1}(z^\alpha-\Delta_h)^{-1}\big\|
\lesssim \tau^2 |z|^{1-\alpha},\quad \forall z \in \Gamma_{\vartheta,\rho}^\tau,
\end{split}\end{equation}
for some $\vartheta \in (\pi/2,\pi)$ and $\rho>0$,
which leads to the following lemma.
\begin{lemma}
The estimate (\ref{cor.3}) holds if and only if $a^{(1)}=-\frac{1}{4}$.
\begin{proof}
For sufficiently small $\tau$, both of the fractional BDF-2 and GNG-2 satisfy (\ref{pre.0}) with $p=2$, indicating that for $|z|\tau<\epsilon$ where $\epsilon$ is chosen sufficiently small, there holds
\begin{equation*}\begin{split}
\big|
z^\alpha-\tau^{-\alpha}\widetilde{\omega}(e^{-z\tau})
\big|\lesssim \tau^2|z|^{2+\alpha}.
\end{split}\end{equation*}
Appealing to the fact
\begin{equation*}\begin{split}
&\quad\tau\Phi_{mi}(e^{-z\tau})
-z^{-1}(z^\alpha-\Delta_h)^{-1}
\\
&=\bigg\{
\frac{2\tau e^{-z\tau}}{1-e^{-z\tau}}
\bigg[
1+a^{(1)}
-\bigg(a^{(1)}+\frac{1}{2}\bigg)
e^{-z\tau}
\bigg]
-z^{-1}
\bigg\}
\big(
\tau^{-\alpha}\widetilde{\omega}(e^{-z\tau})-\Delta_h
\big)^{-1}
\\
&\quad
+z^{-1}
\bigg[
\big(
\tau^{-\alpha}\widetilde{\omega}(e^{-z\tau})
-\Delta_h
\big)^{-1}
-(z^\alpha -\Delta_h)^{-1}
\bigg],
\end{split}\end{equation*}
and using the identity (\ref{bdf.4.2.3}),
one immediately gets
\begin{equation*}\begin{split}
\big|
\tau\Phi_{mi}(e^{-z\tau})
-z^{-1}(z^\alpha-\Delta_h)^{-1}
\big|
&\lesssim
|z|^{-\alpha}
\bigg|
\frac{2\tau e^{-z\tau}}{1-e^{-z\tau}}
\bigg[
1+a^{(1)}
-\bigg(a^{(1)}+\frac{1}{2}\bigg)
e^{-z\tau}
\bigg]
-z^{-1}
\bigg|
\\&\quad
+
\tau^2|z|^{1-\alpha},
\end{split}\end{equation*}
where we have used the estimates (\ref{bdf.4.2}), (\ref{bdf.4.2.1}) and (\ref{pre.4.1}).
Therefore, the estimate (\ref{cor.3}) for sufficiently small $|z|\tau$ holds if and only if
\begin{equation*}\begin{split}
\bigg|
\frac{2\tau e^{-z\tau}}{1-e^{-z\tau}}
\bigg[
1+a^{(1)}
-\bigg(a^{(1)}+\frac{1}{2}\bigg)
e^{-z\tau}
\bigg]
-z^{-1}
\bigg| \lesssim
\tau^2|z|,
\end{split}\end{equation*}
which is equivalent, by resorting to the expansion of the exponential function and some simple calculation, to the fact that $a^{(1)}=-\frac{1}{4}$.
\par
If $|z|\tau > \epsilon$, similar analysis as (ii) or (iii) in Lemma \ref{lem.5} can be developed showing that the estimate (\ref{cor.3}) still holds and the details are omitted here.
The proof of the lemma is completed.
\end{proof}
\end{lemma}
\par
In combination with the relation (\ref{cor.1.1.1}), we obtain $a^{(2)}=\frac{1}{4}$ and the following sharp error estimates.
\begin{theorem}\label{thm.6}
Assume $u_0(x) \in D(\Delta)$ and $f(x,t) \in W^{1,\infty}(0,T;L^2(\Omega))$ and $\int_{0}^{t}(t-s)^{\alpha-1}\|f''(s)\|_{L^2(\Omega)}\mathrm{d}s<\infty$.
Let $W_h^n$ be the solution of (\ref{cor.1}) with $a^{(1)}=-\frac{1}{4}$ and $a^{(2)}=\frac{1}{4}$, and let $U_h^n:=W_h^n+u_h(0)$ be the approximation to $u_h(t_n)$ for $n\geq 1$.
For sufficiently small $\tau>0$, there holds
\begin{equation*}\begin{split}
  \|U_h^n-u_h(t_n)\|_{L^2(\Omega)}
  &\lesssim 
  \tau^2
  \bigg(
  t_n^{\alpha-2}(\|f(0)\|_{L^2(\Omega)}+\|\Delta u_0\|_{L^2(\Omega)})
  +t_n^{\alpha-1}\|f'(0)\|_{L^2(\Omega)}
  \\&\quad
  +\int_{0}^{t_n}(t_n-s)^{\alpha-1}\|f''(s)\|_{L^2(\Omega)}\mathrm{d}s
  \bigg).
\end{split}
\end{equation*}
\begin{proof}
The arguments can be carried out using almost the same strategies as Theorems \ref{thm.2}-\ref{thm.4}, with the estimate (ii) in Lemma \ref{lem.5} replaced by (\ref{cor.3}).
The repetitive details are omitted here.
\end{proof}
\end{theorem}
\begin{remark}
It has been shown in Theorem 3.5 in \cite{jin2013error} that the space semidiscrete solution $u_h(t)$ approximates $u(t)$ with optimal accuracy $O(h^2)$, which in combination with Theorem \ref{thm.6}, leads to the following estimate
\begin{equation*}\begin{split}
\|U_h^n - u(t_n)\|_{L^2(\Omega)} 
\lesssim 
&\tau^2
  \bigg(
  t_n^{\alpha-2}(\|f(0)\|_{L^2(\Omega)}+\|\Delta u_0\|_{L^2(\Omega)})
  +t_n^{\alpha-1}\|f'(0)\|_{L^2(\Omega)}
  \\&
  +\int_{0}^{t_n}(t_n-s)^{\alpha-1}\|f''(s)\|_{L^2(\Omega)}\mathrm{d}s
  \bigg)
+
h^2\|\Delta u_0\|_{L^2(\Omega)}.
\end{split}
\end{equation*}
\end{remark}

\section{Numerical experiments}\label{sec.expr}
In the introduction section, we have mentioned a numerical test with a zero source term showing that the ACN scheme by the fractional BDF-2 or GNG-2 is only $O(\tau^\alpha)$ accuracy, which confirms the results in Theorem \ref{thm.2}.
In this section, more numerical examples will be presented to verify other theoretical results.
We note that only temporal accuracy is reported next and the space mesh size $h$ is taken sufficiently small ($h=10^{-4}$).
Let $\Omega=(0,\pi)$ and $T=1$.
For given time step $\tau$, denoted by $E(t_n,\tau)=\|U_h^n-u(t_n)\|_{L^2(\Omega)}$ the error at $t_n$, or by $E(\tau)=\displaystyle\max_n \|U_h^n-u(t_n)\|_{L^2(\Omega)}$ the error for all $0 \leq n \leq N$.
The convergence order is obtained by
\begin{equation*}
\text{Order}=\log_2 \frac{E(t_n,\tau)}{E(t_n,\tau/2)},\quad \text{or} \quad
\text{Order}=\log_2 \frac{E(\tau)}{E(\tau/2)}.
\end{equation*}
In each table, the theoretical convergence order is offered in parentheses.
\\
\textbf{Example 1.}
In this example, by taking
\begin{equation*}
u(x,t)=\big(E_\alpha(-t^\alpha)-1\big)\sin x \quad \text{and} \quad f(x,t)=-\sin x,
\end{equation*}
where $E_\alpha (t)=\sum_{j=0}^{\infty}\frac{t^j}{\Gamma(\alpha j+1)}$,
one gets $u_0 \equiv 0$.
Numerical results are collected in Table \ref{tab2} for different $\alpha (\alpha=0.2,0.4,0.8)$ for both ACN schemes generated by fractional BDF-2 and GNG-2.
The numerical convergence order indicates the accuracy is of $O(\tau^\alpha)$ which is in line with the theoretical results of Theorem \ref{thm.3}.
\begin{table}[]
\caption{Convergence order at $t_n=0.5$ of ACN scheme for Example 1.}
\label{tab2}
\begin{tabular}{ccccccc}
\hline\noalign{\smallskip}
\multirow{2}{*}{$\alpha$} & \multirow{2}{*}{$\tau$} & \multicolumn{2}{c}{ACN(FBDF-2)} &  & \multicolumn{2}{c}{ACN(GNG-2)} \\ \cline{3-4} \cline{6-7} 
                          &                         & $\|U_h^n-u(t_n)\|_{L^2(\Omega)}$             & Order       &  & $\|U_h^n-u(t_n)\|_{L^2(\Omega)}$            & Order      \\ \hline \noalign{\smallskip}
\multirow{3}{*}{0.2}      & $1/2^7$                 & 2.7962E-01        & (0.20)         &  & 2.7024E-01       & (0.20)         \\
                          & $1/2^8$                 & 2.5066E-01        & 0.16        &  & 2.4202E-01       & 0.16       \\
                          & $1/2^9$                 & 2.2401E-01        & 0.16        &  & 2.1609E-01       & 0.16       \\ \hline \noalign{\smallskip}
\multirow{3}{*}{0.4}      & $1/2^7$                 & 9.5481E-02        & (0.40)          &  & 9.0387E-02       & (0.40)         \\
                          & $1/2^8$                 & 7.3723E-02        & 0.37        &  & 6.9720E-02       & 0.37       \\
                          & $1/2^9$                 & 5.6680E-02        & 0.38        &  & 5.3560E-02       & 0.38       \\ \hline \noalign{\smallskip}
\multirow{3}{*}{0.8}      & $1/2^7$                 & 8.4719E-03        & (0.80)          &  & 8.1956E-03       & (0.80)         \\
                          & $1/2^8$                 & 4.8781E-03        & 0.80        &  & 4.7188E-03       & 0.80       \\
                          & $1/2^9$                 & 2.8059E-03        & 0.80        &  & 2.7142E-03       & 0.80       \\ \hline
\end{tabular}
\end{table}
\\
\textbf{Example 2.}
To further validate Theorem \ref{thm.3}, we choose in this example the source term $f(x,t)=-t\sin x$ and zero initial condition.
The exact solution is unknown and is represented by numerical solutions on fine meshes ($\tau=2^{-12},h=10^{-4}$).
We note that $f(x,t)$ in this case satisfies $f(x,0)=0$, which by Theorem \ref{thm.3}, implies the convergence order is $2$ at any positive time and is $1+\alpha$ if all time levels are considered.
As illustrated in Table \ref{tab3} where errors and convergence orders at time $t_n=0.5$ are reported, one clear observes that for different $\alpha (\alpha=0.2,0.5,0.9)$, the accuracy is of $O(\tau^2)$.
Table \ref{tab4} offers numerical results in the norm $\displaystyle\max_n \|U_h^n - u(t_n)\|_{L^2(\Omega)}$, which confirms the accuracy is of $O(\tau^{1+\alpha})$.
\begin{table}[]
\caption{Convergence order at $t_n=0.5$ of ACN scheme for Example 2.}
\label{tab3}
\begin{tabular}{ccccccc}
\hline\noalign{\smallskip}
\multirow{2}{*}{$\alpha$} & \multirow{2}{*}{$\tau$} & \multicolumn{2}{c}{ACN(FBDF-2)} &  & \multicolumn{2}{c}{ACN(GNG-2)} \\ \cline{3-4} \cline{6-7} 
                          &                         & $\|U_h^n-u(t_n)\|_{L^2(\Omega)}$             & Order       &  & $\|U_h^n-u(t_n)\|_{L^2(\Omega)}$            & Order      \\ \hline\noalign{\smallskip}
\multirow{3}{*}{0.2}      & $1/2^7$                 & 6.3847E-07        & (2.00)          &  & 6.4131E-07       & (2.00)         \\
                          & $1/2^8$                 & 1.5922E-07        & 2.00        &  & 1.5982E-07       & 2.00       \\
                          & $1/2^9$                 & 3.9401E-08        & 2.01        &  & 3.9484E-08       & 2.02       \\ \hline\noalign{\smallskip}
\multirow{3}{*}{0.5}      & $1/2^7$                 & 1.3976E-06        & (2.00)          &  & 1.4668E-06       & (2.00)         \\
                          & $1/2^8$                 & 3.5052E-07        & 2.00        &  & 3.6687E-07       & 2.00       \\
                          & $1/2^9$                 & 8.6878E-08        & 2.01        &  & 9.0731E-08       & 2.02       \\ \hline\noalign{\smallskip}
\multirow{3}{*}{0.9}      & $1/2^7$                 & 1.7237E-06        & (2.00)          &  & 1.5230E-06       & (2.00)         \\
                          & $1/2^8$                 & 4.1513E-07        & 2.05        &  & 3.6628E-07       & 2.06       \\
                          & $1/2^9$                 & 1.0081E-07        & 2.04        &  & 8.8893E-08       & 2.04       \\ \hline
\end{tabular}
\end{table}
\begin{table}[]
\caption{Convergence order of ACN scheme for Example 2 with the maximal norm in time.}
\label{tab4}
\begin{tabular}{ccccccc}
\hline\noalign{\smallskip}
\multirow{2}{*}{$\alpha$} & \multirow{2}{*}{$\tau$} & \multicolumn{2}{c}{ACN(FBDF-2)} &  & \multicolumn{2}{c}{ACN(GNG-2)} \\ \cline{3-4} \cline{6-7} 
                          &                         & $\displaystyle\max_{n}\|U_h^n-u(t_n)\|_{L^2(\Omega)}$             & Order       &  & $\displaystyle\max_{n}\|U_h^n-u(t_n)\|_{L^2(\Omega)}$            & Order      \\ \hline\noalign{\smallskip}
\multirow{3}{*}{0.2}      & $1/2^7$                 & 1.8429E-05        & (1.20)          &  & 1.1902E-05       & (1.20)         \\
                          & $1/2^8$                 & 9.1538E-06        & 1.01        &  & 5.2547E-06       & 1.18       \\
                          & $1/2^9$                 & 4.5100E-06        & 1.02        &  & 2.2767E-06       & 1.21       \\ \hline\noalign{\smallskip}
\multirow{3}{*}{0.5}      & $1/2^7$                 & 4.4319E-05        & (1.50)          &  & 3.1894E-05       & (1.50)         \\
                          & $1/2^8$                 & 1.6699E-05        & 1.41        &  & 1.2143E-05       & 1.39       \\
                          & $1/2^9$                 & 6.1310E-06        & 1.45        &  & 4.4965E-06       & 1.43       \\ \hline\noalign{\smallskip}
\multirow{3}{*}{0.9}      & $1/2^7$                 & 2.1268E-05        & (1.90)          &  & 2.0355E-05       & (1.90)         \\
                          & $1/2^8$                 & 5.8525E-06        & 1.86        &  & 5.5666E-06       & 1.87       \\
                          & $1/2^9$                 & 1.5772E-06        & 1.89        &  & 1.5007E-06       & 1.89       \\ \hline
\end{tabular}
\end{table}
\\
\textbf{Example 3.}
In this example, we consider the modified ACN scheme by taking $u(x,t)$ and $f(x,t)$ as follows
\begin{equation*}
u(x,t)=\big(E_\alpha(-t^\alpha)+t^3\big)\sin x,\quad
f(x,t)=\bigg(\frac{6t^{3-\alpha}}{\Gamma(4-\alpha)}+t^3\bigg)\sin x.
\end{equation*}
We emphasize that although $f(x,t)$ is singular at initial time, it indeed meets the property $f(x,t) \in W^{1,\infty}(0,T;L^2(\Omega))$ and $\int_{0}^{t}(t-s)^{\alpha-1}\|f''(s)\|_{L^2(\Omega)}\mathrm{d}s<\infty$ required by Theorem \ref{thm.6}.
The numerical results are reported in Table \ref{tab5} for different $\alpha (\alpha=0.1,0.5,0.9)$.
Clearly, with the help of corrections at initial two steps, the optimal accuracy $O(\tau^2)$ is arrived at.
\begin{table}[]
\caption{The optimal convergence order at $t_n=0.5$ of modified ACN scheme for Example 3.}
\label{tab5}
\begin{tabular}{ccccccc}
\hline
\multirow{2}{*}{$\alpha$} & \multirow{2}{*}{$\tau$} & \multicolumn{2}{c}{MACN(FBDF-2)} &  & \multicolumn{2}{c}{MACN(GNG-2)} \\ \cline{3-4} \cline{6-7} 
                          &                         & $\|U_h^n-u(t_n)\|_{L^2(\Omega)}$             & Order        &  & $\|U_h^n-u(t_n)\|_{L^2(\Omega)}$            & Order       \\ \hline\noalign{\smallskip}
\multirow{3}{*}{0.1}      & $1/2^7$                 & 7.4149E-07        & (2.00)           &  & 5.9903E-07       & (2.00)          \\
                          & $1/2^8$                 & 1.9277E-07        & 1.94         &  & 1.4581E-07       & 2.04        \\
                          & $1/2^9$                 & 4.4067E-08        & 2.13         &  & 4.1134E-08       & 1.83        \\ \hline\noalign{\smallskip}
\multirow{3}{*}{0.5}      & $1/2^7$                 & 1.0273E-05        & (2.00)           &  & 5.1934E-06       & (2.00)          \\
                          & $1/2^8$                 & 2.6169E-06        & 1.97         &  & 1.3398E-06       & 1.95        \\
                          & $1/2^9$                 & 6.5502E-07        & 2.00         &  & 3.3467E-07       & 2.00        \\ \hline\noalign{\smallskip}
\multirow{3}{*}{0.9}      & $1/2^7$                 & 4.9142E-05        & (2.00)           &  & 4.6469E-05       & (2.00)          \\
                          & $1/2^8$                 & 1.2325E-05        & 2.00         &  & 1.1655E-05       & 2.00        \\
                          & $1/2^9$                 & 3.0805E-06        & 2.00         &  & 2.9127E-06       & 2.00        \\ \hline
\end{tabular}
\end{table}
\section{Conclusion}\label{sec.conc}
In this work, the averaging Crank-Nicolson (ACN) scheme is considered to numerically solving the subdiffusion problems.
Two types of time stepping methods, namely the fractional BDF-2 and the generalized 2nd-order Newton–Gregory formula are adopted to build the ACN scheme.
The generating function involved for such scheme is characterized by its zeros which may be on the unit circle.
By resorting to the residue theorem, sharp error estimates are developed showing that the accuracy of ACN scheme is of $O(\tau^\alpha)$ at any positive time.
To improve the accuracy of the method, corrections are designed at initial two steps which can yield the optimal $O(\tau^2)$ accuracy.
Several numerical examples are conducted to validate all theoretical results obtained in this work.
\begin{acknowledgements}
This work is supported by the Autonomous Region Level High-Level Talent Introduction Research Support Program in 2022 (No. 12000-15042224 to B.Y.)
and National Natural Science Foundation of China (No. 12201322 to B.Y., 12061053 to Y.L.
and 12161063 to H.L.) and Natural Science Foundation of Inner Mongolia (2020MS01003 to Y.L., and 2021MS01018 to H.L.).
\end{acknowledgements}

\bibliographystyle{spmpsci}      
\bibliography{mybibfile}   

\end{document}